\documentclass{amsart}
\usepackage{amsmath}
\usepackage{amsfonts}

\newtheorem{theorem}{Theorem}
\theoremstyle{plain}

\newtheorem{corollary}{Corollary}

\newtheorem{notation}{Notation}
\newtheorem{problem}{Problem}

\newtheorem{remark}{Remark}

\numberwithin{equation}{section}
\input{tcilatex}

\begin{document}
\title[Behavior of Solutions of Nonlinear PDE's and Dynamics of Cancer ]{%
Behavior of Solutions Nonlinear Reaction-Diffusion PDE's Relation to
Dynamics of Propagation of Cancer}
\author{Kamal N. Soltanov}
\address{National Academy of Sciences of Azerbaijan, Baku, Azerbaijan}
\email{sultan\_kamal@hotmail.com }
\urladdr{https://www.researchgate.net/profile/Kamal-Soltanov/research\ }

\begin{abstract}
In this paper, we propose a new mathematical model nonlinear
reaction-diffusion PDE's describing the dynamics of propagation of cancer.
Here the mixed problem for the proposed PDE's is investigated and by
applying obtained results conclusions on the dynamics of propagation of
cancer are drawn. These problems have nonlocal nonlinearity with variable
exponents and possess special properties: these can be to remain either
dissipative all time or become non-dissipative after a finite time. Here the
solvability and behavior of solutions both when problems are yet dissipative
and when become nondissipative are proved. It is shown that if the studied
process gets become nondissipative can have various states, e.g. an infinite
number of different unstable solutions with varying speeds, in addition,
their propagation can become chaotic. The behavior of these solutions is
analyzed in detail and it is explained how space-time chaos can arise.
Investigation of this mathematics model allows explaining the dynamics of
propagation of cancer, which are provided here as conclusions for each case.
\end{abstract}

\maketitle

\section{\protect\bigskip \label{Sec_1}Introduction}

In this article, we introduce and study a mathematical model nonlinear
reaction-diffusion equation describing the dynamics of propagation of
cancer. The studied here partial differential equations (PDE's) are
nonlinear equations of parabolic type with variable exponent nonlocal
nonlinearity that in addition, possess special properties. Here the mixed
problem for this class of nonlinear reaction-diffusion equations is studied.

The theory of nonlinear PDE's of parabolic type is of great interest both in
itself and also as a useful mathematical model for a wide variety of
important problems. This theory is widely used in understudies of various
problems of hydrodynamics, the theory of nonlinear diffusion, and also
understudies many problems in physics, chemistry, biology, etc. For further
information on applications, the reader is referred to such sources as \cite%
{ag,an,ba,cl,cr,du-m,k-m,li,o-m,q,shi,m1,m2,m3,m-p-b,v,z}, etc..

Many biological processes at mathematical modeling are described by the
nonlinear reaction-diffusion equations (reaction-diffusion-convection or
advection). In addition, the mathematical modeling of the dynamics of the
process of diseases due to infection also be described by the
reaction-diffusion equation. It should be noted that for the study of the
dynamics of propagation of cancer usually were used the various problems for
nonlinear reaction-diffusion equations (\cite%
{ag,bek,bel1,bel2,bo,ch,ch-lo,cr,da,de,e,fr1,fr2,lo,k-n,nag,nan,p,sh-fr,t,w}
and the references therein).

It needs to note the nonlinear parabolic PDE's (also reaction-diffusion
equations) comprise infinite-dimensional dynamical systems that exhibit an
amazing spectrum of solution phenomena including traveling waves,
dissipative solitons, spiral waves, target patterns, bifurcation cascades,
chaos, and long-time dynamical configurations of great complexity (\cite%
{an,ba,cl,h,ni,q,shi,m-p-b,t,v,w,z} and the references therein). Especially
interest represents an investigation of the long-time behavior of the
solutions that concerns the dynamic regimes, which the system may settle
into as growing the time. If the system is dissipative, there is typically a
global strange attractor, but mathematically things are especially
challenging when the system is non-dissipative inasmuch as the long-time
dynamics can be much more diverse and complicated than in the dissipative
case. It needs to note that long-time dynamics of solutions of nonlinear
reaction-diffusion equations in the non-dissipative cases have many open
problems. In this article, we obtain new results on the finite and long-time
dynamics of classes of nonlinear reaction-diffusion equations, which provide
additional mathematical details regarding blow-up and chaotic transitions.
Such type questions were analyzed in the literature by many authors (see, 
\cite{ba,bo,cl,du-m,h,li,ni,o-m,q,shi,m-p-b,t,v,z} and the references
therein).

Many authors for the study of the dynamics of propagation of cancer account
that the mathematical model describing this process is a reaction-diffusion
PDE's. Moreover, is assumed that this model will feel the changes of the
state of the cells according to time if the study will be to use the
mathematics model in the form of the mixed problem with the free boundary
that depends on time. According to medical investigations known that genes
can either be activated or suppressed when signals stimulate receptors on
the cell surface and are then transmitted to the nucleus of the cell. The
reception of particular signals can induce a cell to reproduce itself in the
form of identical descendants, that is the so-called clonal expansion or
mitosis, or to die, that is the so-called apoptosis or programmed death
(e.g. \cite%
{ag,al,bek,bel1,bo,cr,da,de,e,fr2,h,lo,ni,t,w,d-n-t,l-t-w-z,p-q-v,k-n} and
the references therein).

We use reaction-diffusion PDE's for the study of the same problem also, but
our approach is different from the models were used in the above-mentioned
works. Here we take account of the mentioned properties by selecting the
coefficients and exponents of nonlinearity as the functions which are
dependents of variables ($t,x$). Since the changes of the state of the cells
pass over time and position, one needs to use such a model that can feel
these changes. According to these properties, it requires to take into
account that the changes of cells are the spatiotemporal process that in
particular, are arose in the appearance of cancer. Consequently, studying
the properties of the solutions of these PDEs helps explain the dynamics of
propagation of the natural processes under discussion in applications. The
here obtained results we used to analyze the dynamics of cancer spread. In
this paper, we provide some explanation of these dynamics and how the
changes can spread over a long-time in various cases.

In what follows we study a mixed problem for an equation that changes from
dissipative to non-dissipative reaction-diffusion PDEs with increasing time.
Here we study also the long-time dynamics of classes of nonlinear
reaction-diffusion equations, which during time pass from dissipative to
non-dissipative equations. Moreover, the corresponding steady-state problem
passes from a uniquely solvable problem to a problem with infinitely many
solutions. It is shown that the trajectories of solutions in the phase space
depend on the choice of starting point on a sphere of the initial values.
Therefore the behavior of solutions to this problem can be more complicated.
Here can arise many variants, e.g. can be the blow-up at finite time, can
exist such absorbing manifold, that the associated dynamics tend to be
chaotic. This class of nonlinear reaction-diffusion equations we study for
analyzing the dynamics of propagation of cancer. \ \ \ \ \ \ \ \ \ \ 

In this article, we investigate a reaction-diffusion problem that over time
is converted to the advection reaction-diffusion process that can have
spatiotemporal chaos. So, we will study the mixed problem for the following
equation as the first step of the approach to the study of the posed
question 
\begin{equation}
\frac{\partial u}{\partial t}=\nabla \cdot \left( a\left( x\right) \nabla
u\right) -b\left( x\right) u-f\left( t,x\right) \left\Vert u\right\Vert
^{p_{0}\left( t\right) }u+g\left( t,x\right) \left\Vert u\right\Vert
^{p_{1}\left( t\right) }u,\ \left( t,x\right) \in R_{+}\times \Omega
\label{1.1}
\end{equation}%
where $\Omega \subset R^{n}$, $n\geq 1$ is the bounded domain with the
sufficiently smooth boundary $\partial \Omega $, \bigskip $u\left(
t,x\right) $ is unknown function and 
\begin{equation*}
\left( i\right) \quad a\left( x\right) ,b\left( x\right) ,f\left( t,x\right)
,g\left( t,x\right) >0,p_{0}\left( t,\left\Vert u_{0}\right\Vert \right)
\geq 0,p_{1}\left( t,\left\Vert u_{0}\right\Vert \right) \geq 0
\end{equation*}%
are known functions.

We presuppose that equation (\ref{1.1}) is the mathematical model that
describes the dynamics of propagation of cancer, on the other words,
describes the corresponding process, which is necessary to investigate. Here
the function $u(t,x)$ denote the density (mass) of cells (as the normal
cells), which are proliferating cells, but these states during time can be
changed according to changes of cells from proliferating to infected (such
as proliferating, quiescent or dead cells). It is necessary to note that $%
\left\Vert u\right\Vert $ denote the density of all cells from the domain $%
\Omega $ and $p_{0}\left( t\right) $ denote the change of the range of
normal cells, and $p_{1}\left( t\right) $ denotes the change of all cells
that do not receive signals from the immune system during at time $t$, i.e.
infected cells.

We should be noted the suggested here mathematical model also is the free
boundary problem but it differs from the usual problems with free boundary
since this model is such as the "two-phase" process where exists the
separating boundary between the normal cells and the destructed cells that
is the free boundary, which generates by the variable exponents $p_{1}\left(
t\right) $ and $p_{0}\left( t\right) $. The equation (\ref{1.1}) assumed if
the function $g\left( t,x\right) >0$ then the change begins to take place.
The transformation of the cells from the normal state to the solid-state is
the known process of vary, therefore the coefficients of the equation (\ref%
{1.1}) can describe by the known functions. So, here for the immune and
"infected" cells, we assume one and the same denotations since these both
are the cells. It needs to note usually authors assume that the process of
cancer already takes place and therefore separately designates each of the
states of cells. According to this approach, the dynamics of the states of
the cells are investigated by the system of equations.

In this article, we use another approach for the study of the same problem.
We will study this process by using such a mathematical model that will feel
the changes of the cells with time and state, which cells receive. Due to
medical investigations is known that genes can either be activated or
suppressed when signals stimulate receptors on the cell surface (and then
transmitted to the nucleus of the cell). These properties we take into
account by selecting the coefficients and exponents as the functions which
are dependent on variables $(t,x)$. The cellular scale refers to the main
(interactive) activities of the cells: activation and proliferation of tumor
cells and competition with immune cells. Since the cells from one state to
another pass gradually, i.e. from the immune state to the activation and
proliferation (interactive) state. The activities (interactiveness) of the
cells after infections mainly are changed in the following sequence:
activation and proliferation of tumor cells and competition with immune
cells. As one can see in the works dedicated to cancer the different states
of cells are usually denoted by different notations (see, e.g. mentioned
above), but the offered here mathematical model feels of these different
states by coefficients and exponents. In the equation (\ref{1.1})$\ a\left(
x\right) $ is the coefficient of diffusion measuring the mobility of any
cells (namely, the immune cells and proliferating cells) and the functions $%
f\left( t,x\right) $ and $g\left( t,x\right) $ are the recruitment rate of
susceptible and degenerated cells, respectively.

So, we will study the following mixed problem with homogeneous boundary
condition: the equation (\ref{1.1}) with the initial and boundary conditions 
\begin{equation}
u\left( 0,x\right) =u_{0}\left( x\right) ,\ u\left( t,x\right) \geq 0,\ 
\label{1.2}
\end{equation}%
\begin{equation}
u\left( t,x\right) \left\vert \ _{\partial \Omega \times R_{+}}=0\right. ,\
p_{0}\left( t\right) \geq 0,\ p_{0}\left( 0\right) =p-1.  \label{1.3}
\end{equation}%
In this problem assumed that the external interferences are excluded.

Let the following conditions (\textit{ii}) 
\begin{equation*}
p_{1}\left( t\right) \geq 0;\ \exists t_{0}>0,\ t>t_{0}\Longrightarrow
p_{1}\left( t\right) >\ p_{0}\left( t\right) ;\exists t_{1}>t_{0},\
t>t_{1}\Longrightarrow \ p_{0}\left( t\right) =0;
\end{equation*}%
hold, where $t_{0}$ defined by equality $\ t_{0}=\inf \left\{ t\in
R_{+}\left\vert \ \frac{p_{1}\left( t\right) }{p_{0}\left( t\right) }%
>1\right. \right\} $. \ Where, in general, $p_{0}\left( t,\tau \right)
\searrow $ and $p_{1}\left( t,\tau \right) \nearrow $, $g\left( t,x\right)
\nearrow $\ if $t,\tau \nearrow $, $b\left( x\right) >$ $g\left( 0,x\right)
\geq 0$. Here $p>0$ denotes the quantity of cells in the examined domain $%
\Omega $.

In this article, we study a new class of the nonlinear reaction-diffusion
PDE's with the nonlocal nonlinearity with the variable exponents. Here we
investigated the solvability and the behavior of the solutions for the
considered problems both when these are yet dissipative and when these get
become non-dissipative. These problems are possessed special properties:
these can be to remain dissipative all time or can be to exists finite time
after that these get become non-dissipative. The long-time behavior of the
solutions things is especially challenging mathematically when the PDE is
non-dissipative inasmuch as these dynamics can be much more diverse and
complicated than in the dissipative case. Here we study the long-time
behavior of the solutions that are necessary to study, as it gives provide
additional details on solutions and especially regarding the chaotic cases.
The dynamics becomes much more complicated in the case PDEs largely due to
the formation of space-time chaotic modes (see e.g. \cite%
{an,ba,li,ni,q,m-a,m-p-b,t,v,w,z}, and the references therein). It should be
noted that there have been many investigations on this and related topics
(see, for example, \ \cite{cl,cr,h,q,m-p-b,t,w,z} and the references
therein). It should be noted that the corresponding steady-state problem of
a non-dissipative nonlinear diffusion-reaction PDE has infinitely many
solutions and because of which there is a bifurcation of solutions,
consequently, appear the chaos (see, e.g. \cite{32}). It is shown that if
the studied problem becomes non-dissipative can arise an infinite number of
different unstable solutions with varying speeds. In this case, due to
bifurcations will generate an infinite number of different states of
spatio-temporal chaos.

We should be noted that in this article our main aim is to study the
behavior of solutions to the posed problem, which is necessary for the study
of the process of the dynamics of propagation of cancer, consequently to
answer when could be the blow-up (or collapse) and what kind can be the
long-time dynamics of propagation of cancer.

This article is constructed in the following way: In Section 2 the
solvability of the problem (\ref{1.1})-(\ref{1.3}) is studied that has 3
subsections where the variants depending on the relation between $p_{0}$ and 
$p_{1}$ are studied. In Section 3 the behavior of solutions of the
considered problem is studied and the blow-up of solutions under finite time
is proved. In Section 4 the long-time dynamics of the solutions of the
examined problem is investigated, where the chaotics of the behavior of
solutions of this problem also is shown. Each section contains the
conclusion about the dynamics of propagation of cancer, where is provided an
explanation: what means to the studied process the obtained in this section
results.

\section{\label{Sec_2}Solvability of problem (\protect\ref{1.1})-(\protect
\ref{1.3})}

Let $\Omega \subset R^{n}$, $n\geq 1$, is bounded domain with sufficiently
smooth boundary $\partial \Omega $ (at least from the Lipschitz class). By $%
W_{0}^{1,2}\left( \Omega \right) \equiv H_{0}^{1}\left( \Omega \right) $ we
denote the Sobolev space and by $L^{r}\left( \Omega \right) $, $r\geq 1$ the
Lebesgue spaces.

Assume the known functions of the problem satisfy the condition (\textit{i})
and are bounded continuous functions, where $Q_{T}=\left( 0,T\right) \times
\Omega $, $T>0$ is a number. It is need noted $\nabla \cdot \left( a\ \nabla
\right) :H_{0}^{1}\left( \Omega \right) \longrightarrow H^{-1}\left( \Omega
\right) $.

Now we cosider the solvability of this problem, which will be analyzed
making use of the general results from \cite{m2,m-p-b,m-a}. We take $%
u_{0}\in B_{r_{0}}^{H_{0}^{1}}\left( 0\right) ,\footnote{%
Where 
\begin{equation*}
B_{r_{0}}^{H_{0}^{1}}\left( 0\right) =\left\{ u\in H_{0}^{1}\left\vert \
\left\Vert \nabla u\right\Vert \leq r_{0}\right. \right\} ,
\end{equation*}%
need to note $\left\Vert \nabla u\right\Vert \equiv $ $\left\Vert
u\right\Vert _{H_{0}^{1}}.$}$ where $r_{0}>0,$ and study the operator $A$
generated by the problem (\ref{1.1})-(\ref{1.3}): it acts, by definition,
from $X:=W^{1,2}\left( 0,T;H^{-1}\left( \Omega \right) \right) \cap $ $%
L^{2}\left( 0,T;H_{0}^{1}\left( \Omega \right) \right) \cap $ $\left\{
u\left( t,x\right) \left\vert \text{ }\ u\left( 0,x\right) =\right.
u_{0}\right\} $ to $Y=L^{2}\left( 0,T;H^{-1}\left( \Omega \right) \right) $.
Consequently, as the solution of the problem (\ref{1.1})-(\ref{1.3}) we
understand follows: a function $u\in X$ called the solution of this problem
if it satisfies this problem in the sense the space $Y=L^{2}\left(
0,T;H^{-1}\left( \Omega \right) \right) $.

According to conditions (\textit{i}) and (\textit{ii}) the equation (\ref%
{1.1}) can changes its character from the dissipative equation to the
non-dissipative, whence ensue that here is possible the following different
cases, which necessary investigate separately: namely, there exist such
times $0\leq t_{0}\leq t_{1}<\infty $ that (1) $p_{0}\left( t\right)
>p_{1}\left( t\right) $, $t\in \left[ 0,t_{0}\right] $; (2) $p_{1}\left(
t\right) \geq p_{0}\left( t\right) $ for $t\in \left[ t_{0},t_{1}\right] $;
(3) $p_{0}\left( t\right) =0$ for $t>t_{1}$. Since the character of the
problem (\ref{1.1})- (\ref{1.3}) also will changes depending on the analysed
case.

\subsection{\label{Ssec_2.1}Existence in the case $p_{0}\left( t\right)
>p_{1}\left( t\right) $}

So, let $p_{0}\left( t\right) >p_{1}\left( t\right) $. In the beginning it
is needed to obtain a priori estimates for the possible solutions, therefore
we consider the following expression \ 
\begin{equation*}
\left\langle \frac{\partial u}{\partial t},u\right\rangle =\left\langle
\nabla \cdot \left( a\nabla u\right) -bu-f\left\Vert u\right\Vert
^{p_{0}\left( t\right) }u+g\left\Vert u\right\Vert ^{p_{1}\left( t\right)
}u,u\right\rangle \Longrightarrow \ 
\end{equation*}%
and in sequal having carried out some calculations we get%
\begin{equation*}
\frac{1}{2}\frac{d}{dt}\left\Vert u\right\Vert _{2}^{2}=-\left\langle
a\nabla u,\nabla u\right\rangle +\left\langle a\frac{\partial u}{\partial
\nu },u\right\rangle \left\vert \ _{\partial \Omega }\right. -\left\langle
bu,u\right\rangle -\left\Vert u\right\Vert ^{p_{0}\left( t\right)
}\left\langle fu,u\right\rangle +
\end{equation*}%
\begin{equation*}
+\left\Vert u\right\Vert ^{p_{1}\left( t\right) }\left\langle
gu,u\right\rangle =-\left\Vert a^{\frac{1}{2}}\nabla u\right\Vert
_{2}^{2}+\left\langle a\frac{\partial u}{\partial \nu },u\right\rangle
\left\vert \ _{\partial \Omega }\right. -
\end{equation*}%
\begin{equation*}
-\left\Vert b^{\frac{1}{2}}u\right\Vert ^{2}-\left\Vert u\right\Vert
^{p_{0}\left( t\right) }\left\Vert f^{\frac{1}{2}}\left( t\right)
u\right\Vert ^{2}+\left\Vert u\right\Vert ^{p_{1}\left( t\right) }\left\Vert
\left( g\right) ^{\frac{1}{2}}\left( t\right) u\right\Vert ^{2}
\end{equation*}%
according of the boundary condition 
\begin{equation*}
\frac{1}{2}\frac{d}{dt}\left\Vert u\right\Vert ^{2}=-\left\Vert a^{\frac{1}{2%
}}\nabla u\right\Vert ^{2}-\left\Vert b^{\frac{1}{2}}u\right\Vert
^{2}-\left\Vert u\right\Vert ^{p_{0}\left( t\right) }\left\Vert f^{\frac{1}{2%
}}u\right\Vert ^{2}+
\end{equation*}%
\begin{equation}
+\left\Vert u\right\Vert ^{p_{1}\left( t\right) }\left\Vert g^{\frac{1}{2}%
}u\right\Vert ^{2},\quad \left\Vert u\right\Vert ^{2}\left( 0\right)
=\left\Vert u_{0}\right\Vert ^{2}  \label{2.1}
\end{equation}%
where $\left\langle \circ ,\circ \right\rangle =\underset{\Omega }{\int }%
\circ \times \circ \ dx$, $\left\Vert \circ \right\Vert $ is the norm of $%
L^{2}\left( \Omega \right) =H\left( \Omega \right) $.

For estimate of the last term in the equation (\ref{2.1}) we consider the
expression $g\left( t,x\right) \cdot \left\Vert u\right\Vert ^{p_{1}\left(
t\right) }$ due of 
\begin{equation*}
\left\Vert u\right\Vert ^{p_{1}\left( t\right) }\left\Vert g^{\frac{1}{2}%
}u\right\Vert ^{2}=\underset{\Omega }{\int }gu^{2}\left\Vert u\right\Vert
^{p_{1}\left( t\right) }dx\Longrightarrow gu^{2}\left\Vert u\right\Vert
^{p_{1}\left( t\right) }.
\end{equation*}

Hence use $\varepsilon -$Young inequality we have \footnote{%
We wish to note when degradation of cells begin $p_{1}\left( t\right) $
remains less than $p_{0}\left( t\right) $ until some time $t_{0}$. By
conditions take place: $p_{1}\left( t\right) \nearrow $ and $p_{0}\left(
t\right) \searrow $ when $t\nearrow $. \ } 
\begin{equation*}
g\left\Vert u\right\Vert ^{p_{1}\left( t\right) }\ \leq \frac{p_{0}-p_{1}}{%
p_{0}}\left( \varepsilon ^{\frac{p_{0}}{p_{0}-p_{1}}}g^{\frac{p_{0}}{%
p_{0}-p_{1}}}f^{-\frac{p_{1}}{p_{0}-p_{1}}}\right) +\frac{p_{1}}{p_{0}}%
\left( \varepsilon ^{\frac{p_{0}}{p_{1}}}f\left\Vert u\right\Vert
^{p_{0}\left( t\right) }\right) ,
\end{equation*}%
and taking $\varepsilon =\left( \frac{p_{1}}{p_{0}}\right) ^{\frac{p_{1}}{%
p_{0}}}$ 
\begin{equation}
g\left\Vert u\right\Vert ^{p_{1}\left( t\right) }\ \leq f\left\Vert
u\right\Vert ^{p_{0}\left( t\right) }+\left[ \left( \frac{p_{1}}{p_{0}}%
\right) ^{\frac{p_{1}}{p_{0}-p_{1}}}-\left( \frac{p_{1}}{p_{0}}\right) ^{%
\frac{p_{0}}{p_{0}-p_{1}}}\right] \left( \frac{g}{f}\right) ^{\frac{p_{0}}{%
p_{0}-p_{1}}}.  \label{2.2}
\end{equation}

If denote by $b_{1}\left( t,x\right) $ the following expression 
\begin{equation*}
b_{1}\left( t,x\right) =b\left( x\right) -\left[ \left( \frac{p_{1}\left(
t\right) }{p_{0}\left( t\right) }\right) ^{\frac{p_{1}\left( t\right) }{%
p_{0}\left( t\right) -p_{1}\left( t\right) }}-\left( \frac{p_{1}\left(
t\right) }{p_{0}\left( t\right) }\right) ^{\frac{p_{0}\left( t\right) }{%
p_{0}\left( t\right) -p_{1}\left( t\right) }}\right] \ \left( \frac{%
g^{p_{0}\left( t\right) }\left( t,x\right) }{f^{p_{1}\left( t\right) }\left(
t,x\right) }\right) ^{\frac{1}{p_{0}\left( t\right) -p_{1}\left( t\right) }}
\end{equation*}%
then from (\ref{2.1}) we derive the inequality 
\begin{equation}
\frac{1}{2}\frac{d}{dt}\left\Vert u\right\Vert ^{2}\leq -\underset{\Omega }{%
\int }a\left( x\right) \left\vert \nabla u\right\vert ^{2}dx-\underset{%
\Omega }{\int }b_{1}\left( t,x\right) u^{2}dx  \label{2.3}
\end{equation}

Hence one need to examine 2 variants: (a) $b_{1}\left( t,x\right) \geq 0$
and (b) $b_{1}\left( t,x\right) <0$. Let takes place (a) then inequation 
\begin{equation}
\underset{\Omega }{\int }a\left( x\right) \left\vert \nabla u\right\vert
^{2}dx+\underset{\Omega }{\int }b_{1}\left( t,x\right) u^{2}dx\geq 0.
\label{2.4}
\end{equation}%
holds and we derive the following problem 
\begin{equation*}
\frac{1}{2}\frac{d}{dt}\left\Vert u\right\Vert ^{2}\leq -\underset{\Omega }{%
\int }a\left( x\right) \left\vert \nabla u\right\vert ^{2}dx-\underset{%
\Omega }{\int }b_{1}\left( t,x\right) u^{2}dx,\ \left\Vert u\right\Vert
^{2}\left( 0\right) =\left\Vert u_{0}\right\Vert ^{2}.
\end{equation*}

Assume there exist such numbers $a_{0},A_{0},b_{0},B_{0}>0$ that $a_{0}\leq
a\left( x\right) \leq A_{0}$ and $b_{0}\leq b\left( x\right) \leq B_{0}$.
Whence if we denote by $b_{1}\left( t\right) =\inf \left\{ b_{1}\left(
t,x\right) \left\vert \ x\in \Omega \right. \right\} $ then we get the
problem 
\begin{equation*}
\frac{1}{2}\frac{d}{dt}\left\Vert u\right\Vert ^{2}\leq -a_{0}\left\Vert
\nabla u\right\Vert ^{2}-b_{1}\left( t\right) \left\Vert u\right\Vert ^{2},\
\left\Vert u\right\Vert ^{2}\left( 0\right) =\left\Vert u_{0}\right\Vert
^{2}.
\end{equation*}

Then one can affirm that solutions will be bounded. Moreover, maybe a
solution will remained stable when $t\nearrow \infty $ if $p_{1}\left(
t\right) $ will remains less than $p_{0}\left( t\right) $ for any $t>0$.

\begin{remark}
Consequently, the investigated process either doesn't make worsen or will
improve and remains in the bounded vicinity of zero when $t\nearrow \infty $%
, i.e. there exists a bounded absorbing subset in the phase space where tend
all trajectories at $t\nearrow \infty $.
\end{remark}

Now let takes place (b) then if we denote by $\widetilde{b}_{1}\left(
t\right) =\sup \left\{ b_{1}\left( t,x\right) \left\vert \ x\in \Omega
\right. \right\} $ then we have the problem 
\begin{equation*}
\frac{1}{2}\frac{d}{dt}\left\Vert u\right\Vert ^{2}\leq -a_{0}\left\Vert
\nabla u\right\Vert ^{2}+\widetilde{b}_{1}\left( t\right) \left\Vert
u\right\Vert ^{2},\ \left\Vert u\right\Vert ^{2}\left( 0\right) =\left\Vert
u_{0}\right\Vert ^{2}.
\end{equation*}

In this case if 
\begin{equation*}
a_{0}\left\Vert \nabla u\right\Vert ^{2}\geq \left\vert \widetilde{b}%
_{1}\left( t\right) \right\vert \left\Vert u\right\Vert ^{2}
\end{equation*}%
then the previous assertion occurs.

Since the case of nonfufilment of (\ref{2.4}) is particular case of the
studied in the subsection 3, here we not will discuss on this.

Thus the following result is proved.

\begin{theorem}
\label{Th_1} Let $p_{1}\left( t\right) <p_{0}\left( t\right) $ for $t\in %
\left[ 0,t_{0}\right) $. Let all above conditions and (\ref{2.4}) are
fulfilled.

Then the problem (\ref{1.1})-(\ref{1.3}) solvable and solutions remain in
the bounded vicinity of the zero for $t\nearrow $.
\end{theorem}

From this theorem and inequalities of such type as (\ref{2.2}) and (\ref{2.4}%
) follows

\begin{corollary}
\label{C_1}If is fulfilled case (1), i.e. $p_{1}\left( t\right) <p_{0}\left(
t\right) $ and $g\left( t,x\right) <f\left( t,x\right) $ for $t\in \left[
0,t_{0}\right) $ then the problem (\ref{1.1})-(\ref{1.3}) solvable and
solutions remain in the bounded vicinity of the zero for $t\nearrow $.
\end{corollary}

The above discussions shows that the solvability of the problem (\ref{1.1})-(%
\ref{1.3}) in this case follows from general results from \cite{m1,m2,m3}.

\subsection{\label{Con_1}Conclusion on the dynamics of cancer}

It is clear that the immune system of an organism act on infected cells to
stop of propagation of the degeneration or reanimate of such cells. If the
immune system is sufficiently strong then it can reanimate certain parts of
degenerated cells or, at least, stop the continuation of the degeneration,
i.e. can change the dynamics of propagation. Consequently, if this is
possible then the exponent $p_{1}\left( t\right) $ can vary and maybe, not
will increase.

If the immune system isn't sufficiently powerful then it can't reanimate
degenerate cells or stop the continuation of the degeneration, i.e. can't
change the dynamics of propagation. If the function $p_{1}\left( t\right) $
continuously increases this shows that this organism can't control the
degeneration and consequently, the propagation of cancer grows. In this
case, is necessary the complements aid of the doctors, i.e. is necessary
control in the form of external actions. In other words, is necessary to use
chemotherapy, radiotherapy, drugs, etc., and also is necessary to reinforce
the immune system of this human.\footnote{%
The control of the dynamics of propagation of cancer will be discussed in
the next paper, where will be provided and some numerical examples.}

\subsection{\label{Ssec_2.2}Analysis of the case $p_{1}\left( t\right)
=p_{0}\left( t\right) $}

According to conditions $p_{0}\left( t\right) \searrow $ and $p_{1}\left(
t\right) \nearrow $ when $t\nearrow $ there exists such time $t_{0}>0$ that $%
p_{1}\left( t_{0}\right) =p_{0}\left( t_{0}\right) $. In this case the
relation between $f\left( t,x\right) $ and $g\left( t,x\right) $ can changed
depending of points of $\Omega $.

It is clear that if $p_{1}\left( t\right) =p_{0}\left( t\right) $ and $%
g\left( t,x\right) \leq f\left( t,x\right) $ for $t>t_{0}$ then the
propagation of the possible solutions will rightly determinable as in the
previous subsection, but the case $g\left( t,x\right) >f\left( t,x\right) $
for $t>t_{0}$ is similar to the case $p_{1}\left( t\right) \geq p_{0}\left(
t\right) $ for $t\geq t_{0}$ therefore, we will be to study it later. But if
the relation between $g\left( t,x\right) $ and $f\left( t,x\right) $ for $%
t>t_{0}$ is undetermined then maybe arise a chaos. Therefore, this case will
be best to investigate together of the case $p_{1}\left( t\right)
>p_{0}\left( t\right) $.

\subsection{\label{Ssec_2.3}Solvability in the case $p_{1}\left( t\right)
>p_{0}\left( t\right) $}

Let the case (3) is fulfilled, i.e. $p_{1}\left( t\right) \geq p_{0}\left(
t\right) $ for $t\geq t_{0}$. If $p_{1}\left( t\right) >p_{0}\left( t\right) 
$ then independent of the relation between $g\left( t,x\right) $ and $%
f\left( t,x\right) $ the behavior of the possible solutions generally
speaking will indeterminable, as their behavior will has vary according of $%
p_{1}\left( t\right) $, of the initial data and the spectrum of the Laplace
operator. It should be noted that according to conditions $f\left(
t,x\right) ,p_{0}\left( t\right) \searrow $ and $g\left( t,x\right)
,p_{1}\left( t\right) \nearrow $ as $t\nearrow $ , therefore acros some time 
$t>t_{1}$ will be: $g\left( t,x\right) >f\left( t,x\right) $ and $%
p_{1}\left( t\right) \gg p_{0}\left( t\right) $.

So, let us $p_{1}\left( t\right) >p_{0}\left( t\right) $ for $t\geq t_{0}$
and again will examine (\ref{2.1}), more exactly \ 
\begin{equation*}
\frac{1}{2}\frac{d}{dt}\left\Vert u\right\Vert ^{2}+\underset{\Omega }{\int }%
\ \left[ a\left( x\right) \left\vert \nabla u\right\vert ^{2}+b\left(
x\right) u^{2}+\left\Vert u\right\Vert ^{p_{0}\left( t\right) }f\left(
t,x\right) u^{2}\right] dx=
\end{equation*}%
\begin{equation}
=\left\Vert u\right\Vert ^{p_{1}\left( t\right) }\underset{\Omega }{\int }\
g\left( t,x\right) u^{2}dx,\quad \left\Vert u\right\Vert ^{2}\left(
t_{0}\right) =\left\Vert u_{t_{0}}\right\Vert ^{2}.  \label{3.1}
\end{equation}

The expression $\left\Vert u\right\Vert ^{p_{0}\left( t\right) }f\left(
t,x\right) $ one can estimate by use of the expression $\left\Vert
u\right\Vert ^{p_{1}\left( t\right) }g\left( t,x\right) $ like in subsection
2.1 
\begin{equation}
f\ \left\Vert u\right\Vert ^{p_{0}\left( t\right) }\leq \frac{p_{1}-p_{0}}{%
p_{1}}\left( \frac{f^{p_{1}}}{g^{p_{0}}}\right) ^{\frac{1}{p_{1}-p_{0}}}+\ 
\frac{p_{0}}{p_{1}}g\ \left\Vert u\right\Vert ^{p_{1}\left( t\right) }.
\label{3.2}
\end{equation}%
Consequently, there exist such functions $b_{1}\left( t,x\right) $ and $%
g_{1}\left( t,x\right) $ that the following inequality 
\begin{equation*}
\frac{1}{2}\frac{d}{dt}\left\Vert u\right\Vert ^{2}+\underset{\Omega }{\int }%
\ \left[ a\left( x\right) \left\vert \nabla u\right\vert ^{2}+b_{1}\left(
t,x\right) u^{2}\right] dx-\left\Vert u\right\Vert ^{p_{1}\left( t\right) }%
\underset{\Omega }{\int }\ g_{1}\left( t,x\right) u^{2}dx\leq 0
\end{equation*}%
holds, where 
\begin{equation*}
b_{1}\left( t,x\right) =b_{1}\left( b,g,f,p_{0},p_{1}\right) ;\quad
g_{1}\left( t,x\right) =g_{1}\left( g,p_{0},p_{1}\right) .
\end{equation*}%
Thus, from (\ref{2.1}) we derive 
\begin{equation}
\frac{1}{2}\frac{d}{dt}\left\Vert u\right\Vert ^{2}+\underset{\Omega }{\int }%
\ \left[ a\left( x\right) \left\vert \nabla u\right\vert ^{2}+b_{1}\left(
t,x\right) u^{2}\right] dx-\left\Vert u\right\Vert ^{p_{1}\left( t\right) }%
\underset{\Omega }{\int }\ g_{1}\left( t,x\right) u^{2}dx\leq 0,
\label{3.3a}
\end{equation}%
\begin{equation}
\quad \left\Vert u\right\Vert ^{2}\left( 0\right) =\left\Vert
u_{0}\right\Vert ^{2}.  \label{3.3}
\end{equation}%
So, from the Cauchy problem (\ref{3.1}) we derive the Cauchy problem (\ref%
{3.3a})-(\ref{3.3}) for the inequality, which permits to investigate of the
behavior of solutions of the problem (\ref{2.1}).

Denote by $g_{10}=\sup \left\{ g_{1}\left( t,x\right) \left\vert \ x\in
\Omega \right. \right\} $ and $b_{10}=\inf \left\{ \left\vert b_{1}\left(
t,x\right) \right\vert \left\vert \ x\in \Omega \right. \right\} ,$ which
where exist due to conditions of this problem. We assume that the Laplace
operator $-\Delta :H_{0}^{1}\left( \Omega \right) \longrightarrow
H^{-1}\left( \Omega \right) $ has only a point spectrum, i.e. $\sigma \left(
-\Delta \right) =\sigma _{p}\left( -\Delta \right) \subset R_{+}$ and denote
by $\lambda _{1}$ the minimal eigenvalue of the Laplace operator $-\Delta $
(we will note that in this case $\left\Vert \nabla u\right\Vert ^{2}\geq
\lambda _{1}\left\Vert u\right\Vert ^{2}$ be valid).

We now investigate the problem (\ref{3.3a})-(\ref{3.3}) for the initial data
satisfying the condition $u_{0}\in S_{r_{0}}^{H_{0}^{1}}\left( 0\right)
\subset H_{0}^{1}\left( \Omega \right) \footnote{%
Where 
\begin{equation*}
S_{r_{0}}^{H_{0}^{1}}\left( 0\right) =\left\{ u\in H_{0}^{1}\left\vert \
\left\Vert \nabla u\right\Vert =r_{0}\right. \right\} .
\end{equation*}%
}$. Under above conditions one can determine the solution of the problem (%
\ref{3.3a})-(\ref{3.3}) formally as 
\begin{equation}
r\left( t\right) ^{2}\leq \exp -\left\{ 2\left( a_{0}\lambda
_{1}+b_{1}\right) t-2g_{10}\underset{0}{\overset{t}{\int }}r\left( \tau
\right) ^{p_{1}\left( \tau \right) }d\tau \right\} r_{0}^{2},  \label{3.4}
\end{equation}%
where $r\left( t\right) =\left\Vert u\right\Vert \left( t\right) $ and $%
r_{0}=\left\Vert u_{0}\right\Vert $. Now if we assume 
\begin{equation*}
g_{10}\left( e^{-\left( a_{0}\lambda _{1}+\overline{b}_{1}\right)
t_{0}}r_{0}\right) ^{p_{1}\left( t_{0}\right) }<a_{0}\lambda
_{1}+b_{10}\left( t_{0}\right) ,
\end{equation*}%
where $\overline{b}_{1}=\sup \left\{ b_{1}\left( t\right) \left\vert \ t\in
\left( 0,t_{0}\right) \right. \right\} $, then from (\ref{3.4}) we get 
\begin{equation*}
\frac{1}{2}\frac{d}{dt}\left\Vert u\right\Vert ^{2}\leq -a_{0}\left\Vert
\nabla u\right\Vert ^{2}-b_{10}\left\Vert u\right\Vert ^{2}+g_{10}\left\Vert
u\right\Vert ^{p_{1}\left( t\right) +2}\leq
\end{equation*}%
\begin{equation*}
-\left( a_{0}\lambda _{1}+b_{10}-g_{10}\left\Vert u\right\Vert ^{p_{1}\left(
t\right) }\right) \left\Vert u\right\Vert ^{2},\quad \left\Vert u\left(
t\right) \right\Vert ^{2}\left\vert _{t=t_{0}}\right. =\left\Vert u\left(
t_{0}\right) \right\Vert ^{2}.
\end{equation*}

Thus we have 
\begin{equation*}
\left\Vert u\right\Vert ^{2}\leq \exp \left\{ -2\left( a_{0}\lambda
_{1}+b_{10}-g_{10}r_{0}^{p_{1}\left( t_{0}\right) }\right) t\right\}
\left\Vert u\left( t_{0}\right) \right\Vert ^{2},
\end{equation*}

Consequently, in this case on the solvability of the problem (\ref{1.1})-(%
\ref{1.3}) one can formulate the following result that follows from the
general results of \cite{m1,m2,m3}.

\begin{theorem}
\label{Th_2}Let all above conditions on the problem (\ref{1.1})-(\ref{1.3})
are fulfilled. Then this problem solvable for any $u_{0}\in
B_{r_{0}}^{H_{0}^{1}}\left( 0\right) \subset H_{0}^{1}\left( \Omega \right) $
if $r_{0}$ satisfies the inequality $g_{10}\left( e^{-\left( a_{0}\lambda
_{1}+\overline{b}_{1}\right) t_{0}}r_{0}\right) ^{p_{1}\left( t_{0}\right)
}<a_{0}\lambda _{1}+b_{10}\left( t_{0}\right) $. Moreover the mapping
(semi-flow) $S\left( t\right) :u_{t_{0}}\longrightarrow u\left( t\right) $
is such that $H$ strongly $S\left( t\right) :\left(
B_{r_{0}}^{H_{0}^{1}}\left( 0\right) \right) \searrow 0$ as $t\nearrow
\infty $.
\end{theorem}

\subsection{\label{Con_3}Conclusion on the dynamics of cancer}

This result shows the rate of the diffusion process in the body can help
immune systems to improve the process of propagation of the beginning to
degenerate cells at a certain moment after the start of the main process.

\section{\label{Sec_3}Behavior of solutions}

Now we will study the case when $p_{0}\left( t\right) \approx 0$ for $t\geq
t_{1}$. Moreover, since the inequality (\ref{3.2}) shows that if $%
p_{1}\left( t\right) >p_{0}\left( t\right) $ then always instead of the
equation from (\ref{3.1}) one can use the inequation (\ref{3.3}), which is
sufficient for the investigation of the behavior of solutions. One can make
use of the formal solution to the problem (\ref{1.1})-(\ref{1.3}) by takes
account of the above assumption onto the known functions. But if the known
functions depend at both of variable then the study of the posed problem
becomes very complicated. Therefore, we will investigate the behavior of
solutions of the problem in a somewhat weak form. As will be seen that such
an approach is enough to receive the necessary information about the
behavior of solutions.

Now we will begin to study the main questions.

\subsection{\label{Ssec_3.1}Blow-up of Solutions}

So, we will begin to investigate the behavior of solutions of posed problem
when $p_{0}\left( t\right) \approx 0$.

Let us $p_{0}\left( t\right) =0$ for $t\geq t_{1}$. We will examine again (%
\ref{2.1}) written in the following form 
\begin{equation*}
\frac{1}{2}\frac{d}{dt}\left\Vert u\right\Vert ^{2}+\underset{\Omega }{\int }%
\ \left[ a\left( x\right) \left\vert \nabla u\right\vert ^{2}+\left( b\left(
x\right) +f\left( t,x\right) \right) u^{2}\right] dx-
\end{equation*}%
\begin{equation}
-\left\Vert u\right\Vert ^{p_{1}\left( t\right) }\underset{\Omega }{\int }\
g\left( t,x\right) u^{2}dx=0,\quad \left\Vert u\right\Vert ^{2}\left(
t_{1}\right) =\left\Vert u_{t_{1}}\right\Vert ^{2}.  \label{4.1}
\end{equation}

According to conditions there exist such constants $%
a_{0},A_{0},b_{0},B_{0},f_{0},F_{0},g_{0},G_{0}>0$ that the following
inequalities 
\begin{equation*}
a_{0}\leq a\left( x\right) \leq A_{0};b_{0}\leq b\left( x\right) \leq
B_{0};f_{0}\leq f\left( t,x\right) \leq F_{0};g_{0}\leq g\left( t,x\right)
\leq G_{0}
\end{equation*}%
hold.

Consequently, we will study the following problems 
\begin{equation}
\frac{1}{2}\frac{d}{dt}\left\Vert u\right\Vert ^{2}+A_{0}\left\Vert \nabla
u\right\Vert ^{2}+\left( B_{0}+F_{0}\right) \left\Vert u\right\Vert
^{2}-g_{0}\left\Vert u\right\Vert ^{p_{1}\left( t\right) +2}\geq 0
\label{4.2}
\end{equation}%
and 
\begin{equation}
\frac{1}{2}\frac{d}{dt}\left\Vert u\right\Vert ^{2}+a_{0}\left\Vert \nabla
u\right\Vert ^{2}+\left( b_{0}+f_{0}\right) \left\Vert u\right\Vert
^{2}-G_{0}\left\Vert u\right\Vert ^{p_{1}\left( t\right) +2}\leq 0
\label{4.3}
\end{equation}%
with the initial condition $\left\Vert u\right\Vert ^{2}\left( t_{1}\right)
=\left\Vert u_{t_{1}}\right\Vert ^{2}$, which can present the behavior of
the possible solutions of the problem (\ref{4.1}). More exactly, the
trajectory of the possible solutions of the problem (\ref{4.1}) in phase
space go-between of the boundary layers descibed by solutions of the above
problems.

So, from problem (\ref{4.2}) and (\ref{4.3})\ we derive the following
problems 
\begin{equation}
\frac{1}{2}\frac{d}{dt}\left\Vert u\right\Vert ^{2}+A_{0}\left\Vert \nabla
u\right\Vert ^{2}+B_{1}\left\Vert u\right\Vert ^{2}-g_{0}\left\Vert
u\right\Vert ^{p_{1}\left( t\right) +2}\geq 0  \label{4.4}
\end{equation}%
\begin{equation}
\frac{1}{2}\frac{d}{dt}\left\Vert u\right\Vert ^{2}+a_{0}\left\Vert \nabla
u\right\Vert ^{2}+b_{1}\left\Vert u\right\Vert ^{2}-G_{0}\left\Vert
u\right\Vert ^{p_{1}\left( t\right) +2}\leq 0,  \label{4.5}
\end{equation}%
with the initial condition $\left\Vert u\right\Vert ^{2}\left( t_{1}\right)
=\left\Vert u_{t_{1}}\right\Vert ^{2}$, where $B_{1}=B_{0}+F_{0}$ and $%
b_{1}=b_{0}+f_{0}$.

\begin{remark}
The above inequalities show that it is necessary used the derivative of the
function with the variable exponent, therefore here we deduce it. Consider
the function $y\left( t\right) ^{-q\left( t\right) }$ then 
\begin{equation*}
\frac{d}{dt}\left( y\left( t\right) ^{-q\left( t\right) }\right) =\frac{d}{dt%
}\exp \left\{ -q\left( t\right) \ln y\left( t\right) \right\} =-\exp \left\{
-q\left( t\right) \ln y\left( t\right) \right\} \frac{d}{dt}\left( q\left(
t\right) \ln y\left( t\right) \right) =
\end{equation*}%
\begin{equation*}
=-\exp \left\{ -q\left( t\right) \ln y\left( t\right) \right\} \left[ q^{%
{\acute{}}%
}\left( t\right) \ln y\left( t\right) +\frac{q\left( t\right) }{y\left(
t\right) }y^{%
{\acute{}}%
}\left( t\right) \right] =
\end{equation*}%
\begin{equation}
=y\left( t\right) ^{-q\left( t\right) }\left[ q^{%
{\acute{}}%
}\left( t\right) \ln y\left( t\right) +\frac{q\left( t\right) }{y\left(
t\right) }y^{%
{\acute{}}%
}\left( t\right) \right] .  \label{4.6}
\end{equation}
\end{remark}

Consequently, for investigate of the inequality (\ref{4.5}) we will rewrite
it in the following form 
\begin{equation*}
\frac{1}{2}\frac{d}{dt}\left\Vert u\right\Vert ^{2}+\left\Vert u\right\Vert
^{2}\left( \frac{p_{1}%
{\acute{}}%
\left( t\right) }{2p\left( t\right) }\ln \left\Vert u\right\Vert ^{2}\right)
\leq -a_{0}\left\Vert \nabla u\right\Vert ^{2}-b_{1}\left\Vert u\right\Vert
^{2}+
\end{equation*}%
\begin{equation}
+G_{0}\left\Vert u\right\Vert ^{p_{1}\left( t\right) +2}+\left\Vert
u\right\Vert ^{2}\left( \frac{p_{1}%
{\acute{}}%
\left( t\right) }{2p_{1}\left( t\right) }\ln \left\Vert u\right\Vert
^{2}\right) .  \label{4.7}
\end{equation}

\begin{remark}
Consider the function $h\left( z\right) =z^{\frac{1}{2}}-\ln z$ for $z>0$,
then we have $\frac{d}{dt}h\left( z\right) =\frac{1}{2}z^{-\frac{1}{2}}-%
\frac{1}{z}=\frac{1}{z}\left( \frac{1}{2}z^{\frac{1}{2}}-1\right) $. Whence
for $\frac{d}{dt}h\left( z\right) =0$ we get $z^{\frac{1}{2}}-2=0$ or $z=4$
is the minimum of function $h\left( z\right) $, i.e. $\min \left\{ h\left(
z\right) \left\vert \ z>0\right. \right\} =2-\ln 4>0$. Consequently, $%
h\left( z\right) >0$ for $z>0$ (see, e.g. the proof of Lemma 9 of \cite{m-u}%
).
\end{remark}

Thus we get 
\begin{equation}
\left\Vert u\right\Vert ^{2}\left( \frac{p_{1}%
{\acute{}}%
\left( t\right) }{2p_{1}\left( t\right) }\ln \left\Vert u\right\Vert
^{2}\right) \leq \left\Vert u\right\Vert ^{3}\leq \left( b_{0}+f_{0}\right)
\left\Vert u\right\Vert ^{2}+c\left\Vert u\right\Vert ^{p_{1}\left( t\right)
+2}  \label{4.8}
\end{equation}%
as $p_{1}\left( t\right) \geq 1$, where $c=c\left( b_{0},f_{0},p_{1}\right) $%
. More exactly, use the $\varepsilon -$Young inequality and selecting 
\begin{equation*}
\varepsilon =\frac{p_{1}^{%
{\acute{}}%
}}{2p_{1}}\left[ \left( \frac{p_{1}}{p_{1}-1}\right) \left(
b_{0}+f_{0}\right) \right] ^{-\frac{p_{1}-1}{p_{1}}}=\frac{p_{1}^{%
{\acute{}}%
}}{2p_{1}}\left[ \left( \frac{p_{1}}{p_{1}-1}\right) b_{1}\right] ^{-\frac{%
p_{1}-1}{p_{1}}}
\end{equation*}
we obtain $c=\frac{\varepsilon ^{p_{1}}}{p_{1}}$.

Using the inequality (\ref{4.8}) in (\ref{4.7}) we derive the following
problem 
\begin{equation}
\frac{1}{2}\frac{d}{dt}\left\Vert u\right\Vert ^{2}+\left\Vert u\right\Vert
^{2}\left( \frac{p_{1}%
{\acute{}}%
\left( t\right) }{2p_{1}\left( t\right) }\ln \left\Vert u\right\Vert
^{2}\right) \leq -a_{0}\left\Vert \nabla u\right\Vert ^{2}+\left(
G_{0}+c\right) \left\Vert u\right\Vert ^{p_{1}\left( t\right) +2}
\label{4.9a}
\end{equation}%
with the initial condition 
\begin{equation*}
\left\Vert u\right\Vert ^{2}\left( t_{1}\right) =S\left( t_{1}\right)
\left\Vert u_{0}\right\Vert ^{2}=
\end{equation*}%
\begin{equation}
=\exp \left\{ -\left[ a_{0}\lambda _{1}t_{1}+\overset{t_{0}}{\underset{0}{%
\int }}b_{10}\left( s\right) ds+\overset{t_{1}}{\underset{t_{0}}{\int }}%
\left( b_{1}\left( s\right) -g_{10}\left( s\right) r_{0}^{p_{1}\left(
s\right) }\right) ds\right] \right\} \left\Vert u_{0}\right\Vert ^{2}
\label{4.9b}
\end{equation}%
where $u_{0}\in H_{0}^{1}$ is a given initial function, e.g. one can to
choose it from the ball $B_{\widetilde{r}}^{H_{0}^{1}}\left( 0\right)
\subset H_{0}^{1}$ with the radius $\widetilde{r}>0$. So, we study this
problem for the initial function $u_{0}$ belonging to the $B_{\widetilde{r}%
}^{H_{0}^{1}}\left( 0\right) $. In the previous section was proved that
problem has solution in appropriate space for $t\in \left( 0,t_{1}\right] $,
consequently we can affirm the value of $\left\Vert u\right\Vert ^{2}\left(
t\right) $ is defined on $t_{1}$ as in (\ref{4.9b}) with the semi-flow $%
S\left( t\right) $.

Thus if one denote $\left\Vert u\right\Vert ^{2}\left( t\right) =y\left(
t\right) $, $p\left( t\right) =\frac{p_{1}\left( t\right) }{2}$ and account $%
\left\Vert \nabla u\right\Vert ^{2}\geq \lambda _{1}\left\Vert u\right\Vert
^{2}$ then (\ref{4.9a})-(\ref{4.9b}) we derive the following Cauchy problem 
\begin{equation*}
\frac{1}{2}\frac{dy}{dt}+y\left( \frac{p%
{\acute{}}%
\left( t\right) }{p_{1}\left( t\right) }\ln y\right) \leq -a_{0}\lambda
_{1}y+G_{0}y^{p\left( t\right) +1},\quad y\left( t_{1}\right) =y_{t_{1}}.
\end{equation*}%
here is denoted $p\left( t\right) =\frac{p_{1}\left( t\right) }{2}$. Whence
follows 
\begin{equation*}
p\left( t\right) y^{-p\left( t\right) -1}\frac{dy}{dt}+y^{-p\left( t\right)
}\left( p%
{\acute{}}%
\left( t\right) \ln y\right) \leq -a_{0}\lambda _{1}p_{1}\left( t\right)
y^{-p\left( t\right) }+G_{0}p_{1}\left( t\right)
\end{equation*}%
and consequently, we have the inequality 
\begin{equation*}
-\frac{d\left( y^{-p\left( t\right) }\right) }{dt}\leq -a_{0}\lambda
_{1}p_{1}\left( t\right) y^{-p\left( t\right) }+G_{0}p_{1}\left( t\right)
\end{equation*}%
denoted by $z\left( t\right) =y^{-p\left( t\right) }$ we get 
\begin{equation*}
\frac{dz}{dt}\geq a_{0}\lambda _{1}p_{1}\left( t\right) z-G_{0}p_{1}\left(
t\right) ,\quad z\left( t_{1}\right) =y\left( t_{1}\right) ^{-p\left(
t_{1}\right) }.
\end{equation*}

Hence follows 
\begin{equation*}
z\left( t\right) \geq e^{a_{0}\lambda _{1}\underset{t_{1}}{\overset{t}{\int }%
}p_{1}\left( s\right) ds}z\left( t_{1}\right) -\overset{t}{\underset{t_{1}}{%
\int }}e^{a_{0}\lambda _{1}\underset{\tau }{\overset{t}{\int }}p_{1}\left(
s\right) ds}G_{0}p_{1}\left( \tau \right) d\tau .
\end{equation*}%
\qquad \qquad\ Assume $\frac{d}{dt}P_{1}\left( t\right) =\ p_{1}\left(
t\right) $ then we obtain 
\begin{equation*}
z\left( t\right) \geq e^{a_{0}\lambda _{1}\left( P_{1}\left( t\right)
-P_{1}\left( t_{1}\right) \right) }z\left( t_{1}\right) +G_{0}\underset{t_{1}%
}{\overset{t}{\int }}e^{a_{0}\lambda _{1}\left( P_{1}\left( t\right)
-P_{1}\left( \tau \right) \right) }p_{1}\left( \tau \right) d\tau G_{0}
\end{equation*}%
Consequently, we arrive 
\begin{equation*}
z\left( t\right) \geq e^{a_{0}\lambda _{1}\left( P_{1}\left( t\right)
-P_{1}\left( t_{1}\right) \right) }z\left( t_{1}\right) +\left[
1-e^{a_{0}\lambda _{1}\left( P_{1}\left( t\right) -P_{1}\left( t_{1}\right)
\right) }\right] \frac{G_{0}}{a_{0}\lambda _{1}}=
\end{equation*}%
\begin{equation*}
=\frac{G_{0}}{a_{0}\lambda _{1}}-e^{a_{0}\lambda _{1}\left( P_{1}\left(
t\right) -P_{1}\left( t_{1}\right) \right) }\left[ \frac{G_{0}}{a_{0}\lambda
_{1}}-z\left( t_{1}\right) \right] .
\end{equation*}

Takes account that $z\left( t\right) =y^{-p\left( t\right) }$ we get 
\begin{equation*}
y^{-p\left( t\right) }\geq \frac{G_{0}}{a_{0}\lambda _{1}}-e^{a_{0}\lambda
_{1}\left( P_{1}\left( t\right) -P_{1}\left( t_{1}\right) \right) }\left[ 
\frac{G_{0}}{a_{0}\lambda _{1}}-y\left( t_{1}\right) ^{-p\left( t_{1}\right)
}\right]
\end{equation*}%
or 
\begin{equation*}
y^{p\left( t\right) }\leq \frac{e^{-a_{0}\lambda _{1}\left( P_{1}\left(
t\right) -P_{1}\left( t_{1}\right) \right) }y\left( t_{1}\right) ^{p\left(
t_{1}\right) }}{\left\{ 1-\left[ 1-e^{-a_{0}\lambda _{1}\left( P_{1}\left(
t\right) -P_{1}\left( t_{1}\right) \right) }\right] \frac{G_{0}y\left(
t_{1}\right) ^{p\left( t_{1}\right) }}{a_{0}\lambda _{1}}\right\} }=
\end{equation*}%
\begin{equation*}
=\frac{e^{-a_{0}\lambda _{1}\left( P_{1}\left( t\right) -P_{1}\left(
t_{1}\right) \right) }a_{0}\lambda _{1}y\left( t_{1}\right) ^{p\left(
t_{1}\right) }}{\left\{ a_{0}\lambda _{1}-\left[ 1-e^{-a_{0}\lambda
_{1}\left( P_{1}\left( t\right) -P_{1}\left( t_{1}\right) \right) }\right]
G_{0}y\left( t_{1}\right) ^{p\left( t_{1}\right) }\right\} }.
\end{equation*}%
Therefore one need to investigate the following equation 
\begin{equation*}
a_{0}\lambda _{1}-G_{0}y\left( t_{1}\right) ^{p\left( t_{1}\right)
}+e^{-a_{0}\lambda _{1}\left( P_{1}\left( t\right) -P_{1}\left( t_{1}\right)
\right) }G_{0}y\left( t_{1}\right) ^{p\left( t_{1}\right) }=0
\end{equation*}%
or 
\begin{equation*}
P_{1}\left( t\right) =P_{1}\left( t_{1}\right) -\frac{1}{a_{0}\lambda _{1}}%
\ln \left( 1-\frac{a_{0}\lambda _{1}}{G_{0}y\left( t_{1}\right) ^{p\left(
t_{1}\right) }}\right)
\end{equation*}%
since $P_{1}\left( t\right) $ is the increasing function its inverse
function $\left( P_{1}\left( t\right) \right) ^{-1}$ exists therefore we can
derive the upper bound time of the blow-up.

If one lead the above-mentioned calculations for the Cauchy problem posed
for inequality (\ref{4.4}) then we derive the following inequation 
\begin{equation*}
y^{p\left( t\right) }\geq \frac{e^{-\left( A_{0}\lambda _{1}+\widetilde{B}%
_{0}\right) \left( P_{1}\left( t\right) -P_{1}\left( t_{1}\right) \right)
}\left( A_{0}\lambda _{1}+\widetilde{B}_{0}\right) y\left( t_{1}\right)
^{p\left( t_{1}\right) }}{\left\{ \left( A_{0}\lambda _{1}+\widetilde{B}%
_{0}\right) -\left[ 1-e^{-\left( A_{0}\lambda _{1}+\widetilde{B}_{0}\right)
\left( P_{1}\left( t\right) -P_{1}\left( t_{1}\right) \right) }\right]
G_{0}y\left( t_{1}\right) ^{p\left( t_{1}\right) }\right\} }.
\end{equation*}%
Consequently, we arrive to a result similar to the obtained above result.%
\footnote{%
Consequently, if we wish to obtain such a result as in the above section
then necessary to select the initial function in the appropriate way as in
the above section.} Whence by similar way as of the previous case we can
determine the lower bound time of the blow-up. Therefore the time of blow-up 
$t_{\func{col}}$ of the main problem can be found between these times that
is a finite since $p$ is finite number, and $p_{1}\left( t\right) $ is a
continuous function and satisfies condition $p_{1}\left( t\right) \leq p-1$.

Accordingly, we have proved the following result.

\begin{theorem}
\label{Th_4}Let us the initial function $u_{0}\in H_{0}^{1}$ is such that
the condition 
\begin{equation*}
\left\Vert a^{\frac{1}{2}}\nabla u\left( t_{1}\right) \right\Vert
^{2}-\left\Vert g^{\frac{1}{2}}\left( t_{1}\right) u\left( t_{1}\right)
\right\Vert ^{2}\left\Vert u\left( t_{1}\right) \right\Vert ^{p_{1}\left(
t_{1}\right) }<0
\end{equation*}%
is fulfilled. Then each solution $u\in X$ of the problem (\ref{1.1}) -\ (\ref%
{1.3}) has blow-up in $H$ a finite time.
\end{theorem}

\begin{notation}
Whence one can derive condition on $u\left( 0\right) $ using of the adduced
above definition of $u\left( t_{1}\right) $.
\end{notation}

\subsection{\label{Con_2}Conclusion on the dynamics of cancer}

The result of Theorem 3 shows the following case is possible: There exist
such initial dates that if the process of the dynamics of propagation of
cancer starts with such initial values then beginning on some time almost
everywhere in the part of the body where take place the process will non
remain of the normal cells. Consequently, one can count that the treatment
already is impossible. Whence follows that is necessary exterior
interference, moreover during and before of the critic moment. In the case
when all conditions of this theorem are fulfilled then according to its
result one can define the approximate time when will be already late. We
should note in this case also is necessary to take account of Conclusion \ref%
{Con_1}.

\section{\label{Sec_4}Long-time Behavior of Solutions}

Consider the equation (\ref{1.1}) in the case when $p_{0}\left( t\right) =0$

\begin{equation*}
\frac{\partial u}{\partial t}=\nabla \cdot \left( a\left( x\right) \nabla
u\right) -\left[ b\left( x\right) +f\left( t,x\right) \right] u+g\left(
t,x\right) \left\Vert u\right\Vert ^{p_{1}\left( t\right) }u.\ 
\end{equation*}

For study of the behavior of solutions we will investigate the following
problem 
\begin{equation}
\frac{\partial u}{\partial t}=\widehat{a}\Delta u-\left[ \widehat{b}+%
\widehat{f}\left( t\right) \right] u+\widehat{g}\left( t\right) \left\Vert
u\right\Vert ^{p_{1}\left( t\right) }u,\quad u\left( t_{1}\right) =u_{t_{1}},
\label{5.1}
\end{equation}%
where instead of the functions $a,b,f$ and $g$ set their mean values that
denoted by $\widehat{a},\widehat{b},\widehat{f}$ and $\widehat{g}$, e.g. $%
\widehat{a}=\frac{1}{mes\Omega }\underset{\Omega }{\int }a\left( x\right) dx$%
, in addition we denote $\widehat{b}+\widehat{f}\left( t\right) $ by $%
\widehat{b}\left( t\right) $. We would be noted the function $p_{1}\left(
t\right) $ increases at $t\nearrow $ and converges to the value $p-1$ that
was defined in the introduction.

It should be noted according to the results of \cite{cor,m1} on the
differential operators of elliptic type, a study of the behavior of
solutions of the problem for the equation (\ref{5.1}) will give the
sufficient informations about the behavior of solutions and also for the
main problem.

So, for the study of the behavior of solutions to the problem (\ref{5.1})
consider the following problem 
\begin{equation}
\frac{1}{2}\frac{d}{dt}\left\Vert u\right\Vert ^{2}+\widehat{a}\left\Vert
\nabla u\right\Vert ^{2}+\widehat{b}\left( t\right) \left\Vert u\right\Vert
^{2}-\widehat{g}\left( t\right) \left\Vert u\right\Vert ^{p_{1}\left(
t\right) +2}=0,  \label{5.1a}
\end{equation}%
\begin{equation}
\left\Vert u\left( t_{1}\right) \right\Vert ^{2}=\left\Vert
u_{t_{1}}\right\Vert ^{2},\quad \text{where }\left\Vert u\left( t_{1}\right)
\right\Vert ^{2}=S\left( t_{1}\right) \left\Vert u_{0}\right\Vert ^{2}
\label{5.1b}
\end{equation}%
the semi-flow $S\left( t\right) $ defined in (\ref{4.9b}).

The examined here problem one can compare the studied in the article \cite%
{m-p-b}. It isn't difficult to see that these problems are different
essentially since in the examined problem, unlike the above-mentioned,
coefficients and exponents depend on the independent variables, which create
additional difficulties.

We will use the norm $\left\Vert \nabla u\right\Vert $ for the norm of the
space $H_{0}^{1}\left( \Omega \right) $ due to the equivalence $\left\Vert
u\right\Vert _{H_{0}^{1}\left( \Omega \right) }\equiv \left\Vert \nabla
u\right\Vert $ for $\forall u\in H_{0}^{1}\left( \Omega \right) $. Here for
simplicity we assume the domain $\Omega $ in a geometric sense is such that
Laplace operator $-\Delta :H_{0}^{1}\left( \Omega \right) \longrightarrow
H^{-1}\left( \Omega \right) $ has only a point spectrum, i.e. $\sigma \left(
-\Delta \right) \equiv \sigma _{P}\left( -\Delta \right) \subset \left(
0,\infty \right) $ and we will denote of these eigenvalues of the Laplace
operator $-\Delta :H_{0}^{1}\left( \Omega \right) \longrightarrow
H^{-1}\left( \Omega \right) $ by 
\begin{equation*}
\lambda _{k},\quad k=1,2,...;\quad \sigma _{P}\left( -\Delta \right) \equiv
\left\{ \lambda _{j}\left\vert \ j\in 
\mathbb{N}
\right. \right\} ,
\end{equation*}%
and the corresponding eigenfunctions by 
\begin{equation*}
w_{k},\quad k=1,2,...;\quad \left\{ w_{k}\right\} _{k=1}^{\infty }\subset
H_{0}^{1}\left( \Omega \right) .
\end{equation*}

Now we will study the inverse mapping of the operator $\widehat{A}$
generated by this problem that generated similarly as the inverse mapping of
the operator $A$ of the problem (\ref{5.1a}) -\ (\ref{5.1b}) (i.e. the
mapping $\widehat{S\left( t\right) }$ is defined as the mapping $S\left(
t\right) $ in Section \ref{Sec_2}).

Moreover, we assume that the eigenfunctions and adjoint eigenfunctions are
total (complete) in the space $H_{0}^{1}\left( \Omega \right) $ and also in
the dual space, respectively; in addition, we assume without loss of
generality that they generate an orthogonal basis in these spaces,
respectively.

We introduce the following denotation 
\begin{equation*}
\mathrm{\inf }\left\{ \lambda _{k}\in \sigma _{P}\left( -\Delta \right)
\left\vert \ \widehat{a}\lambda _{k}+\widehat{b}\left( t_{1}\right) >%
\widehat{g}\left( t\right) r\left( t_{1}\right) ^{p_{1}\left( t_{1}\right)
},\right. k=1,2,...\right\} =\lambda _{k_{0}}
\end{equation*}%
$\ \ $ and use the usual representation of the space $H_{0}^{1}\left( \Omega
\right) $ (\cite{an,m1,m-p-b,q,z}) in the form $H_{0}^{1}\left( \Omega
\right) \equiv H_{k_{0}}\oplus H_{-k_{0}}$, where the subspace $%
H_{k_{0}}\subset H_{0}^{1}\left( \Omega \right) $ is the span over $\left\{
w_{k}\right\} _{k=1}^{k_{0}-1}$ and has dimension $\dim H_{k_{0}}=k_{0}-1$
and $H_{-k_{0}}$ is a subspace of co$\mathrm{\dim }H_{-k_{0}}=k_{0}-1.$ Let $%
Q_{k_{0}}$ and $P_{k_{0}}$ are the projections: $P_{k_{0}}:H_{0}^{1}\left(
\Omega \right) \longrightarrow H_{k_{0}}\subset H_{0}^{1}\left( \Omega
\right) $ and $Q_{k_{0}}:H_{0}^{1}\left( \Omega \right) \longrightarrow
H_{-k_{0}}\subset H_{0}^{1}\left( \Omega \right) $, giving rise to the
splitting $u\equiv Q_{k_{0}}u+P_{k_{0}}u$. (This is well known decomposition
of the Hilbert space, see e.g., \cite{m1,m-p-b,q,ni}, etc.).

So, it is easy to see that $-\Delta :H_{k_{0}}\longrightarrow H_{k_{0}}^{-}$
and $-\Delta :H_{-k_{0}}\longrightarrow H_{-k_{0}}^{-},$ where the subspaces 
$H_{k_{0}}^{-},H_{-k_{0}}^{-}$ possess bi-orthogonal bases (see, for
instance, \cite{q,m1,m-p-b,v,z}, etc.), consequently, due to the evident
commutativity of operators $P_{k_{0}}$ and $Q_{k_{0}}$ with the Laplacian $%
\Delta $ \ in $H_{0}^{1}\left( \Omega \right) ,$ one can rewrite the problem
as 
\begin{equation}
\frac{\partial \ }{\partial t}P_{k_{0}}u-\widehat{a}\Delta P_{k_{0}}u+%
\widehat{b}\left( t\right) P_{k_{0}}u-\widehat{g}\left( t\right) \left\Vert
u\right\Vert ^{p_{1}\left( t\right) }P_{k_{0}}u=0,  \label{5.2}
\end{equation}%
\begin{equation}
P_{k_{0}}u\left( t_{1},x\right) =P_{k_{0}}u_{t_{1}}\left( x\right) \in
H_{k_{0}}\subset H_{0}^{1}\left( \Omega \right) ,  \label{5.3}
\end{equation}%
\begin{equation}
\frac{\partial }{\partial t}Q_{k_{0}}u-\widehat{a}\Delta Q_{k_{0}}u+\widehat{%
b}\left( t\right) Q_{k_{0}}u-\widehat{g}\left( t\right) \left\Vert
u\right\Vert ^{p_{1}\left( t\right) }Q_{k_{0}}u=0,  \label{5.4}
\end{equation}%
\begin{equation}
Q_{k_{0}}u\left( t_{1},x\right) =Q_{k_{0}}u_{t_{1}}\left( x\right) \in
H_{-k_{0}}\subset H_{0}^{1}\left( \Omega \right) .  \label{5.5}
\end{equation}

As our aim is the investigation of the behavior of solutions of the problem
under the condition $u_{t_{1}}\in B_{r\left( t_{1}\right) }^{H_{0}^{1}\left(
\Omega \right) }\left( 0\right) $ then it is enough to assume that $\widehat{%
a}\lambda _{k_{0}-1}+\widehat{b}\left( t_{1}\right) <\widehat{g}\left(
t_{1}\right) r\left( t_{1}\right) ^{p_{1}\left( t_{1}\right) }<\widehat{a}%
\lambda _{k_{0}}+\widehat{b}\left( t_{1}\right) $. Then for the problem\ (%
\ref{5.1a}) - \ (\ref{5.1b}) we obtain 
\begin{equation*}
0=\frac{1}{2}\frac{d}{dt}\left\Vert P_{k_{0}}u\right\Vert _{2}^{2}\left(
t\right) +
\end{equation*}%
\begin{equation}
+\left[ \widehat{a}\left\Vert \nabla P_{k_{0}}u\right\Vert ^{2}+\widehat{b}%
\left( t\right) \left\Vert P_{k_{0}}u\right\Vert ^{2}-\widehat{g}\left(
t\right) \left\Vert u\right\Vert ^{p_{1}\left( t\right) }\left\Vert
P_{k_{0}}u\right\Vert ^{2}\right] \left( t\right) ,  \label{5.6a}
\end{equation}%
\begin{equation}
\left\langle P_{k_{0}}u,P_{k_{0}}u\right\rangle \left\vert \
_{t=t_{1}}\right. =\left\Vert P_{k_{0}}u\right\Vert ^{2}\left( t_{1}\right)
=\left\Vert P_{k_{0}}u_{t_{1}}\right\Vert ^{2}.  \label{5.6b}
\end{equation}%
Whence, it follows that the solution $\left\Vert u\right\Vert \left(
t\right) $ and exponent $p_{1}\left( t\right) $ grow and for some $%
t_{2}>t_{1}$ for $t\in \left[ t_{1},t_{2}\right) $ we have 
\begin{equation*}
\widehat{g}\left( t\right) \left\Vert u\right\Vert ^{p_{1}\left( t\right)
}\left( t\right) \leq \widehat{g}\left( t_{1}\right) r\left( t_{1}\right)
^{p_{1}\left( t\right) }+\delta <\widehat{a}\lambda _{k_{0}}+\widehat{b}%
\left( t\right) ,
\end{equation*}%
for some $\delta >0$, by virtue of the inequality: $\widehat{g}\left(
t_{1}\right) r\left( t_{1}\right) ^{p_{1}\left( t_{1}\right) }-\left( 
\widehat{a}\lambda _{k_{0}-1}+\widehat{b}\left( t_{1}\right) \right) >0$.
Indeed if $\left\Vert u_{t_{1}}\right\Vert =r\left( t_{1}\right) $, then we
have from \ (\ref{5.6a}) - \ (\ref{5.6b}) that 
\begin{equation*}
\frac{d}{dt}\left\Vert P_{k_{0}}u\right\Vert ^{2}\left( t\right) +2\left( 
\widehat{a}\lambda _{k_{0}-1}+\widehat{b}\left( t\right) -\widehat{g}\left(
t\right) r\left( t_{1}\right) ^{p_{1}\left( t\right) }\right) \left\Vert
P_{k_{0}}u\right\Vert ^{2}\left( t\right) \geq
\end{equation*}%
\begin{equation*}
\frac{d}{dt}\left\Vert P_{k_{0}}u\right\Vert _{2}^{2}\left( t\right) +2(%
\widehat{a}\lambda _{k_{0}-1}+\widehat{b}\left( t\right) -\widehat{g}\left(
t\right) r\left( t_{1}\right) ^{p_{1}\left( t\right) })\left\Vert
P_{k_{0}}u\right\Vert ^{2}\left( t\right) \geq 0
\end{equation*}%
and we get the inequality 
\begin{equation}
\left\Vert P_{k_{0}}u\right\Vert _{2}^{2}\left( t\right) \geq \exp \left\{ -2%
\underset{t_{1}}{\overset{t}{\int }}\left( \widehat{a}\lambda _{k_{0}-1}+%
\widehat{b}\left( s\right) -\widehat{g}\left( s\right)
r_{t_{1}}^{p_{1}\left( s\right) }\right) ds\right\} \left\Vert
P_{k_{0}}u_{t_{1}}\right\Vert _{2}^{2}.  \label{5.7}
\end{equation}

Consequently if $\left\Vert Q_{k_{0}}u\right\Vert ^{2}\left( t\right) \leq
\epsilon <\delta <r_{0}$ for some enough small $\epsilon >0$ and $t\in \left[
t_{1},t_{2}\right) ,$ then the solution of problem (\ref{5.6a}) - \ (\ref%
{5.6b}) exists and is an exponentially increasing function, by virtue of the
inequality: $\widehat{g}\left( t_{1}\right) r\left( t_{1}\right)
^{p_{1}\left( t_{1}\right) }-\left( \widehat{a}\lambda _{k_{0}}+\widehat{b}%
\left( t_{1}\right) \right) <0$.

Unlike of the above problem, for the problem\ (\ref{5.4}) -\ (\ref{5.5}) we
obtain 
\begin{equation*}
\frac{d}{dt}\left\Vert Q_{k_{0}}u\right\Vert ^{2}\left( t\right) +2\left( 
\widehat{a}\lambda _{k_{0}}+\widehat{b}\left( t\right) -\widehat{g}\left(
t\right) r\left( t_{1}\right) ^{p_{1}\left( t\right) }\right) \left\Vert
Q_{k_{0}}u\right\Vert ^{2}\left( t\right) \leq 0
\end{equation*}%
\begin{equation*}
\left\langle Q_{k_{0}}u,Q_{k_{0}}u\right\rangle \left\vert \
_{t=t_{1}}\right. =\left\Vert Q_{k_{0}}u\right\Vert _{2}^{2}\left(
t_{1}\right) =\left\Vert Q_{k_{0}}u_{t_{1}}\right\Vert _{2}^{2}
\end{equation*}%
whence we get the inequality 
\begin{equation}
\left\Vert Q_{k_{0}}u\right\Vert _{2}^{2}\left( t\right) \leq \exp \left\{ -2%
\underset{t_{1}}{\overset{t}{\int }}\left( \widehat{a}\lambda _{k_{0}-1}+%
\widehat{b}\left( s\right) -\widehat{g}\left( s\right)
r_{t_{1}}^{p_{1}\left( s\right) }\right) ds\right\} \left\Vert
Q_{k_{0}}u_{t_{1}}\right\Vert _{2}^{2}.  \label{5.8}
\end{equation}%
The inequality (\ref{5.8}) shows that the solution of the problem\ (\ref{5.4}%
) -\ (\ref{5.5}) exists and is an exponentially decreasing function, since $%
\widehat{g}\left( t_{1}\right) r\left( t_{1}\right) ^{p_{1}\left(
t_{1}\right) }-\widehat{a}\lambda _{k_{0}}+\widehat{b}\left( t_{1}\right) <0$%
. Consequently for $\left\Vert P_{k_{0}}u_{t_{1}}\right\Vert +\left\Vert
Q_{k_{0}}u_{t_{1}}\right\Vert =r\left( t_{1}\right) $ if $\left\Vert
Q_{k_{0}}u_{t_{1}}\right\Vert <\left\Vert P_{k_{0}}u_{t_{1}}\right\Vert $
and $\left\Vert Q_{k_{0}}u_{t_{1}}\right\Vert $ is enough small, then the
solution $\left\Vert u\right\Vert \left( t\right) $ exists and is an
increasing function up to some time.

Now we will assume the system of eigenfunctions $\left\{ w_{k}\right\}
_{k=1}^{\infty }\subset H_{0}^{1}\left( \Omega \right) $ are an orthonormal
basis of this space that can give possible for detail to investigate of the
behavior of solutions. In this case each function $u\left( t,x\right) $ $\in 
$ $L^{2}\left( \left( 0,T\right) ;H_{0}^{1}\left( \Omega \right) \right) $
has the representation $u\left( t,x\right) =\underset{k=1}{\overset{\infty }{%
\sum }}u_{k}\left( t\right) w_{k}\left( x\right) $ and consequently, study
of the problem (\ref{5.1a}) - (\ref{5.1b}) is equivalent to studying the
system of equations 
\begin{equation*}
\frac{1}{2}\frac{d}{dt}\left\vert u_{k}\left( t\right) \right\vert
^{2}+\left( \widehat{a}\lambda _{k}+\widehat{b}\left( t\right) \right)
\left\vert u_{k}\left( t\right) \right\vert ^{2}-
\end{equation*}%
\begin{equation}
-\widehat{g}\left( t\right) \left( \underset{i=1}{\overset{\infty }{\sum }}%
\left\vert u_{k}\left( t\right) \right\vert ^{2}\right) ^{\frac{p_{1}\left(
t\right) }{2}}\left\vert u_{k}\left( t\right) \right\vert ^{2}=0\ ,\quad
k=1,2,...\quad \   \label{5.9}
\end{equation}%
with the initial conditions 
\begin{equation}
\left\vert u_{k}\left( t_{1}\right) \right\vert ^{2}=\left\vert
u_{t_{1}k}\right\vert ^{2},\quad k=1,2,....  \label{5.9a}
\end{equation}

From mentioned above reasoning follows that for $u\left( t_{1},x\right) \in
B_{r\left( t_{1}\right) }^{H_{0}^{1}\left( \Omega \right) }\left( 0\right) $
and $\left\Vert u_{t_{1}}\right\Vert \leq r\left( t_{1}\right) $ will be
fulfill 
\begin{equation*}
\left\Vert u\right\Vert ^{2}\left( t\right) \equiv \left( \underset{i=1}{%
\overset{\infty }{\sum }}\left\vert u_{i}\left( t\right) \right\vert
^{2}\right) \leq r\left( t_{1}\right) ^{2}+\varepsilon ,
\end{equation*}%
where $t=t\left( \varepsilon ,r\left( t_{1}\right) ,p_{1}\left( t_{1}\right)
\right) >t_{1}$ and $\delta >0$ are sufficiently small. It is known that in
this case $\left\vert u_{k}\left( t\right) \right\vert ^{2}$ increases for $%
k=1,2,...,\widetilde{k}_{0}<$ $k_{0}$ and decreases for $k=k_{0},k_{0}+1,...$
depending on the relationship between $\left\Vert u_{t_{1}}\right\Vert
^{p_{1}\left( t_{1}\right) }$ and $\lambda _{k}$.

Consequently, it is need to investigate the behavior of $\left\vert
u_{k}\left( t\right) \right\vert $ for each $k=1,2,...$. Assume $%
u_{t_{1}}\in H_{0}^{1}\left( \Omega \right) $ and $\left\Vert
u_{t_{1}}\right\Vert \equiv r\left( t_{1}\right) >0$. Let us list all of
possible cases: 1) $\widehat{g}\left( t_{1}\right) r\left( t_{1}\right)
^{p_{1}\left( t_{1}\right) }<\widehat{a}\lambda _{1}+\widehat{b}\left(
t_{1}\right) $; 2) 
\begin{equation*}
\exists \lambda _{k_{0}}:\widehat{a}\lambda _{k_{0}-1}+\widehat{b}\left(
t_{1}\right) <\widehat{g}\left( t_{1}\right) r\left( t_{1}\right)
^{p_{1}\left( t_{1}\right) }<\widehat{a}\lambda _{k_{0}}+\widehat{b}\left(
t_{1}\right)
\end{equation*}%
and 3) $\exists \lambda _{k_{0}}:\widehat{g}\left( t_{1}\right) r\left(
t_{1}\right) ^{p_{1}\left( t_{1}\right) }=\widehat{a}\lambda _{k_{0}}+%
\widehat{b}\left( t_{1}\right) $. Case 1) was already investigated,
therefore we will consider here only cases 2) and 3).

Consider either the case 2) or 3), i.e. 
\begin{equation*}
\exists \lambda _{k_{0}}:\widehat{a}\lambda _{k_{0}-1}+\widehat{b}\left(
t_{1}\right) <\widehat{g}\left( t_{1}\right) r\left( t_{1}\right)
^{p_{1}\left( t_{1}\right) }<\widehat{a}\lambda _{k_{0}}+\widehat{b}\left(
t_{1}\right)
\end{equation*}%
and $\exists \lambda _{k_{0}}:\widehat{g}\left( t_{1}\right) r\left(
t_{1}\right) ^{p_{1}\left( t_{1}\right) }=\widehat{a}\lambda _{k_{0}}+%
\widehat{b}\left( t_{1}\right) .$

So, we will investigate the Cauchy problem for the following system of
equations

\begin{equation*}
\frac{1}{2}\frac{d}{dt}\left\vert u_{k}\left( t\right) \right\vert
^{2}+\left( \widehat{a}\lambda _{k}+\widehat{b}\left( t\right) \right)
\left\vert u_{k}\left( t\right) \right\vert ^{2}-\widehat{g}\left( t\right)
r\left( t\right) ^{p_{1}\left( t\right) }\left\vert u_{k}\left( t\right)
\right\vert ^{2}=
\end{equation*}%
\begin{equation}
=\frac{1}{2}\frac{d}{dt}\left\vert u_{k}\left( t\right) \right\vert
^{2}+\left( \widehat{a}\lambda _{k}+\widehat{b}\left( t\right) -\widehat{g}%
\left( t\right) r\left( t\right) ^{p_{1}\left( t\right) }\right) \left\vert
u_{k}\left( t\right) \right\vert ^{2}=0,  \label{5.10}
\end{equation}%
with the initial conditions 
\begin{equation*}
\left\vert u_{k}\left( t_{1}\right) \right\vert ^{2}=\left\vert
u_{t_{1}k}\right\vert ^{2},k=1,2,...,
\end{equation*}%
where $u\left( t,x\right) \equiv \underset{k=1}{\overset{\infty }{\sum }}%
u_{k}\left( t\right) w_{k}\left( x\right) $. It is easy to see that this
system of equations are such that the cases $k\geq k_{0}$, $k\leq k_{0}-1$
and $k=k_{0}$ is necessary to study separately after which to investigate
the fact that they will be to give when these are together.

Formally, we can determine the solution of each equation from (\ref{5.10})
to be 
\begin{equation}
\left\vert u_{k}\left( t\right) \right\vert ^{2}=\exp \left\{ -2\overset{t}{%
\underset{t_{1}}{\int }}\left( \widehat{a}\lambda _{k}+\widehat{b}\left(
s\right) -\widehat{g}\left( s\right) r\left( s\right) ^{p_{1}\left( s\right)
}\right) ds\right\} \left\vert u_{t_{1}k}\right\vert ^{2}\ .  \label{5.11}
\end{equation}%
So, let us $k\leq k_{0}-1$ then in the case 2), it follows from (\ref{5.11})
that 
\begin{equation*}
\left\vert u_{k}\left( t\right) \right\vert ^{2}=\quad
\end{equation*}%
\begin{equation*}
=\exp \left\{ -2\left( \overset{t}{\underset{t_{1}}{\int }}\left( \widehat{a}%
\lambda _{k}+\widehat{b}\left( s\right) -\widehat{g}\left( s\right) r\left(
s\right) ^{p_{1}\left( s\right) }\right) ds\right) \right\} \left\vert
u_{t_{1}k}\right\vert ^{2}\geq
\end{equation*}%
\begin{equation*}
\geq \exp \left\{ 2\left( \widehat{g}\left( t_{1}\right) r\left(
t_{1}\right) ^{p_{1}\left( t_{1}\right) }-\widehat{a}\lambda _{k}-\widehat{b}%
\left( t_{1}\right) \right) \left( t-t_{1}\right) \right\} \left\vert
u_{t_{1}k}\right\vert ^{2},
\end{equation*}%
for$\ 1\leq k\leq k_{0}-1$ and some $t>t_{1}$ due to $\widehat{g}\left(
t_{1}\right) r\left( t_{1}\right) ^{p_{1}\left( t_{1}\right) }>\widehat{a}%
\lambda _{k}+\widehat{b}\left( t_{1}\right) $ and $p_{1}\left( t\right)
\nearrow $. Consequently, $\left\vert u_{k}\left( t\right) \right\vert $
increases for each $k:$ $1\leq k\leq k_{0}-1$ that lead to the increase of $%
r\left( t\right) $ as long as $\left\Vert P_{k_{t_{1}}}u_{t_{1}}\right\Vert $
is enough greater than $\left\Vert Q_{k_{t_{1}}}u_{t_{1}}\right\Vert $.

For the case 3), i.e. when there exists $k=k_{0}$ such that $\widehat{g}%
\left( t_{1}\right) r\left( t_{1}\right) ^{p_{1}\left( t_{1}\right) }=%
\widehat{a}\lambda _{k_{0}}+\widehat{b}\left( t_{1}\right) $ then one has 
\begin{equation*}
\left\vert u_{k_{0}}\left( t\right) \right\vert ^{2}=\exp \left\{ -2\overset{%
t}{\underset{t_{1}}{\int }}\left( \widehat{a}\lambda _{k_{0}}+\widehat{b}%
\left( s\right) -\widehat{g}\left( s\right) r\left( s\right) ^{p_{1}\left(
s\right) }\right) ds\right\} \left\vert u_{t_{1}k_{0}}\right\vert ^{2}
\end{equation*}%
according to \ (\ref{5.11}). If we denote here the function $\rho \left(
t\right) =\widehat{a}\lambda _{k_{0}}+\widehat{b}\left( t\right) -\widehat{g}%
\left( t\right) r\left( t\right) ^{p_{1}\left( t\right) }$\ then $\rho
\left( t_{1}\right) =0$, but in general the vary of function $\rho \left(
t\right) $ is not known. Consequently, it is impossible to obtain a result
about variation of a solution to the equation (\ref{5.10}), as the behavior
of the function $\rho \left( t\right) $ is not known. As one can see in the
sequel, the behavior of $\rho \left( t\right) $ depends on the geometrical
location of the initial data $u_{t}$ on spheres $S_{r}^{H_{0}^{1}}\left(
0\right) $, $0<r\leq r\left( t_{1}\right) $.

As the functional $\rho \left( t\right) $ depends on $r\left( t\right)
^{p_{1}\left( t\right) }$, $r^{2}\left( t\right) =\left\Vert u\right\Vert
^{2}\left( t\right) $ and $\lambda _{k}$ it is clear that in order to study
the behavior of the functional $\rho \left( t\right) $ one should
investigate both of functionals $\left\Vert P_{k_{0}}u\right\Vert $ and $%
\left\Vert Q_{k_{0}}u\right\Vert $. Clearly, in the case 2) the functional $%
\left\Vert P_{k_{0}}u\right\Vert $ increases, and $\left\Vert
Q_{k_{0}}u\right\Vert $ decreases at $t>t_{1}$ at least near zero according
to \ (\ref{5.6a}) and \ (\ref{5.8}). Using the orthogonal splitting $%
u=P_{k_{0}}u+Q_{k_{0}}u,$ we also find that 
\begin{equation*}
\left\Vert u\right\Vert ^{2}=\left\Vert P_{k_{0}}u\right\Vert
^{2}+\left\Vert Q_{k_{0}}u\right\Vert ^{2}.
\end{equation*}%
Whence follows, the behavior of the function $\left\Vert u\right\Vert
^{2}\left( t\right) $ depends on the relationship between the values $%
\left\Vert P_{k_{0}}u_{t_{1}}\right\Vert $ and $\left\Vert
Q_{k_{0}}u_{t_{1}}\right\Vert .$ Let%
\begin{equation*}
\left\Vert u_{t_{1}}\right\Vert \equiv r\left( t_{1}\right) ;\ \widehat{a}%
\lambda _{k_{0}-1}+\widehat{b}\left( t_{1}\right) <\widehat{g}\left(
t_{1}\right) r\left( t_{1}\right) ^{p_{1}\left( t_{1}\right) }<\widehat{a}%
\lambda _{k_{0}}+\widehat{b}\left( t_{1}\right)
\end{equation*}%
and consider above equality, i.e. 
\begin{equation}
\left\Vert u\right\Vert ^{2}\left( t\right) =\left\Vert
P_{k_{0}}u\right\Vert ^{2}\left( t\right) +\left\Vert Q_{k_{0}}u\right\Vert
^{2}\left( t\right) =\underset{k\leq k_{0}}{\sum }\left\vert u_{k}\left(
t\right) \right\vert ^{2}+\underset{k>k_{0}}{\sum }\left\vert u_{k}\left(
t\right) \right\vert ^{2}.  \label{5.12}
\end{equation}%
It is necessary to study the following cases: \ \textit{a}) 
\begin{equation*}
u_{t_{1}}\equiv \underset{k\leq k_{0}-1}{\sum }u_{t_{1}k}w_{k}\in
P_{k_{0}}\left( H_{0}^{1}\left( \Omega \right) \right) \equiv H_{k_{0}};
\end{equation*}%
\ \textit{b}) 
\begin{equation*}
u_{t_{1}}\equiv \underset{k\geq k_{0}}{\sum }u_{t_{1}k}w_{k}\in
Q_{k_{0}}\left( H_{0}^{1}\left( \Omega \right) \right) \equiv H_{-k_{0}}
\end{equation*}%
and \textit{c}) $\ \ u_{t_{1}}\equiv \underset{k\geq 1}{\sum }%
u_{t_{1}k}w_{k} $ , when \ \textit{c}$_{1}$) $\left\Vert
Q_{k_{0}}u_{t_{1}}\right\Vert <\left\Vert P_{k_{0}}u_{t_{1}}\right\Vert $
and \ \textit{c}$_{2}$) $\left\Vert Q_{k_{0}}u_{t_{1}}\right\Vert \geq
\left\Vert P_{k_{0}}u_{t_{1}}\right\Vert $, \ separately.

In the case \textit{a}) we have 
\begin{equation*}
\left\vert u_{k}\left( t\right) \right\vert ^{2}=\exp \left\{ -2\left( 
\overset{t}{\underset{t_{1}}{\int }}\left( \widehat{a}\lambda _{k}+\widehat{b%
}\left( s\right) -\widehat{g}\left( s\right) r\left( s\right) ^{p_{1}\left(
s\right) }\right) ds\right) \right\} \left\vert u_{t_{1}k}\right\vert ^{2}
\end{equation*}%
for any $k=1,...,k_{0}-1$, thus , $u_{k}\left( t\right) =0$ for $k\geq k_{0}$
since $\left\Vert u_{t_{1}}\right\Vert ^{2}=r\left( t_{1}\right) ^{2}\equiv 
\underset{k\leq k_{0}-1}{\sum }\left( u_{t_{1}k}\right) ^{2}$ and $r\left(
t\right) ^{2}\equiv \underset{k\leq k_{0}-1}{\sum }\left\vert u_{k}\left(
t\right) \right\vert ^{2}$. In this case $\widehat{a}\lambda _{k_{0}}+%
\widehat{b}\left( t\right) -\widehat{g}\left( t\right) r\left( t\right)
^{p_{1}\left( t\right) }<0$ as $\widehat{g}\left( t\right) r\left( t\right)
^{p_{1}\left( t\right) }>\widehat{a}\lambda _{k_{0}}+\widehat{b}\left(
t\right) $ for each $k=1,...,k_{0}-1$ and $p_{1}\left( t\right) $, and also $%
r\left( t\right) $ increase as $t\uparrow \infty $.

In the case b) we have $\left\Vert u_{t_{1}}\right\Vert ^{2}=r\left(
t_{1}\right) ^{2}\equiv \underset{k\geq k_{0}}{\sum }\left(
u_{t_{1}k}\right) ^{2}$, i.e. $u_{t_{1}}\in Q_{k_{0}}\left( \boldsymbol{H}%
\right) \equiv H_{-k_{0}}$. Consequently, 
\begin{equation*}
\left\vert u_{k}\left( t\right) \right\vert ^{2}=\exp \left\{ -2\left( 
\overset{t}{\underset{t_{1}}{\int }}\left( \widehat{a}\lambda _{k}+\widehat{b%
}\left( s\right) -\widehat{g}\left( s\right) r\left( s\right) ^{p_{1}\left(
s\right) }\right) ds\right) \right\} \left\vert u_{t_{1}k}\right\vert ^{2}
\end{equation*}%
for any $k\geq k_{0},$ giving rise to $u_{k}\left( t\right) \searrow 0$ for $%
k=1,2,...,k_{0}-1,$and $r\left( t\right) ^{p_{1}\left( t\right) }$ decreases
although $p_{1}\left( t\right) \nearrow $\ as $t\uparrow \infty $ for each $%
k\geq k_{0}$ due to assumption that $\left( \widehat{a}\lambda _{k_{0}}+%
\widehat{b}\left( s\right) -\widehat{g}\left( s\right) r\left( s\right)
^{p_{1}\left( s\right) }\right) >0$ as $\widehat{g}\left( t\right) r\left(
t\right) ^{p_{1}\left( t\right) }<\widehat{a}\lambda _{k_{0}}+\widehat{b}%
\left( t\right) $. Thus, from continuity we obtain that the inequality $%
\left\Vert u\right\Vert \left( t\right) =r\left( t\right) <r\left(
t_{1}\right) $ and $k\geq k_{0}$ is fulfilled for all $t>t_{1}$ and for any
solution to the problem 
\begin{eqnarray*}
\frac{1}{2}\frac{d}{dt}\left\vert u_{k}\left( t\right) \right\vert
^{2}+\left( \widehat{a}\lambda _{k_{0}}+\widehat{b}\left( t\right) -\widehat{%
g}\left( t\right) r\left( t\right) ^{p_{1}\left( t\right) }\right)
\left\vert u_{k}\left( t\right) \right\vert ^{2} &=&0,\  \\
\left\vert u_{k}\left( t_{1}\right) \right\vert ^{2} &=&\left\vert
u_{t_{1}k}\right\vert ^{2}.
\end{eqnarray*}%
It should be noted above in this section essentially was discussed the case
c) that is general case, since here taken into account both side of the
following inequality 
\begin{equation*}
\widehat{a}\lambda _{k_{0}}+\widehat{b}\left( t_{1}\right) >\widehat{g}%
\left( t_{1}\right) r\left( t_{1}\right) ^{p_{1}\left( t_{1}\right) }>%
\widehat{a}\lambda _{k_{0}-1}+\widehat{b}\left( t_{1}\right) .
\end{equation*}%
Unlike of above cases here it is needed to consider the space decomposition $%
\boldsymbol{H}\equiv H_{k_{0}}\oplus H_{-k_{0}}\equiv H_{0}^{1}\left( \Omega
\right) $ \ \ ($P_{k_{0}}\left( \boldsymbol{H}\right) \equiv H_{k_{0}}$, $%
Q_{k_{0}}\left( \boldsymbol{H}\right) \equiv H_{-k_{0}}$), which is
necessary for analyzes the intrinsic behavior of the corresponding
solutions. Thus, using the decomposition of the space $\boldsymbol{H}\equiv
H_{0}^{1}\left( \Omega \right) $\ for the $r^{2}\left( t_{1}\right)
=\left\Vert u_{t_{1}}^{-}\right\Vert _{\boldsymbol{H}}^{2}$ we have the
representation $r^{2}\left( t_{1}\right) \equiv \left\Vert
u_{t_{1}k_{0}}^{-}\right\Vert _{H_{k_{0}}}^{2}+\left\Vert
u_{t_{1}k_{0}}^{+}\right\Vert _{H_{-k_{0}}}^{2},$ and in the vector form 
\begin{equation*}
\boldsymbol{H}\equiv \left\{ u=\left( u_{k_{0}}^{+},u_{k_{0}}^{-}\right) :\
u_{k_{0}}^{+}\in H_{k_{0}},u_{k_{0}}^{-}\in H_{-k_{0}}\right\} .
\end{equation*}%
Then the following proposition follows directly from the above discussion.

\begin{theorem}
\label{Th_3}Under the above conditions if $\left\Vert
u_{t_{1}k_{0}}^{-}\right\Vert _{H_{-k_{0}}}^{2}>\left\Vert
u_{t_{1}k_{0}}^{+}\right\Vert _{H_{k_{0}}}^{2}$ and if the rate of decrease
of the norm $\left\Vert u_{k_{0}}^{-}\left( t\right) \right\Vert
_{H_{-k_{0}}}^{2}$ is greater than the rate of \ decrease of the$\ $norm $\
\left\Vert u_{k_{0}}^{+}\left( t\right) \right\Vert _{H_{k_{0}}}^{2},$ then
there exists $\widetilde{t}>t_{1}$ such that $\left\vert u_{k}\left(
t\right) \right\vert $ decreases for $t\geq \widetilde{t}$. On the other
hand, if $\left\Vert u_{t_{1}k_{0}}^{-}\right\Vert
_{H_{-k_{0}}}^{2}<\left\Vert u_{t_{1}k_{0}}^{+}\right\Vert _{H_{k_{0}}}^{2}$
and the rate of increase rate of the norm $\left\Vert u_{k_{0}}^{+}\left(
t\right) \right\Vert _{H_{k_{0}}}^{2}$ is greater than the rate of decrease
of $\ \left\Vert u_{k_{0}}^{-}\left( t\right) \right\Vert _{H_{-k_{0}}}^{2}$
, then there exists $\widehat{t}>t_{1}$ such that $\left\vert u_{k}\left(
t\right) \right\vert $ increases for $t\geq \widehat{t}$. Here the increases
of the function $p_{1}\left( t\right) $ loosen only somewhat of the
corresponding behaviors of the solutions.
\end{theorem}

\begin{remark}
This shows how is need to control in order to change of the behavior of the
dynamics of cancer, i.e. to which coefficients or exponents to control is
necessary.
\end{remark}

Now we will proceed to investigate in detail the behavior of solutions
beginning with the time $t_{1}$ when $p_{0}\left( t\right) =0$ (i.e. when
the degregated cells remain only). If one selects an initial function $%
u_{0}\in \boldsymbol{H}$ then into time $t_{1}$ we have the function 
\begin{equation*}
\left\Vert u_{t_{1}}\right\Vert ^{2}=r^{2}\left( t_{1}\right) =
\end{equation*}%
\begin{equation*}
=.\exp \left\{ -\left[ a_{0}\lambda _{1}t_{1}+\overset{t_{0}}{\underset{0}{%
\int }}b_{10}\left( s\right) ds+\overset{t_{1}}{\underset{t_{0}}{\int }}%
\left( b_{1}\left( s\right) -g_{10}\left( s\right) \left( r_{0}\right)
^{p_{1}\left( s\right) }\right) ds\right] \right\} \left\Vert
u_{0}\right\Vert ^{2}
\end{equation*}%
according of the expression (\ref{4.9b}). So, it isn't difficult to see that
arbitrarily chosen initial function $u_{0}$ in the moment $t_{1}$ will be $%
u_{t_{1}}$ with the norm determined by the number $r\left( t_{1}\right) $
that will satisfies the inequalities\ 
\begin{equation*}
\widehat{a}\lambda _{k_{0}}+\widehat{b}\left( t_{1}\right) \geq \widehat{g}%
\left( t_{1}\right) r\left( t_{1}\right) ^{p_{1}\left( t_{1}\right) }>%
\widehat{a}\lambda _{k_{0}-1}+\widehat{b}\left( t_{1}\right)
\end{equation*}%
for some $\lambda _{k_{0}}$.\ As $\boldsymbol{H}=H_{k_{0}}\oplus H_{-k_{0}}$
therefore each $u\in \boldsymbol{H}$ has the representation by decomposition 
$u\left( t\right) =u^{+}\left( t\right) +u^{-}\left( t\right) $, where $%
u^{+}\left( t\right) \in H_{k_{0}}$ and $u^{-}\left( t\right) \in H_{-k_{0}}$
for any $t>0$. As show the above discussions it suffices to study of the
behavior of solutions beginning with the cases that depend of the relations
between of $\left\Vert u_{t_{1}}^{-}\right\Vert _{H_{-k_{0}}}$ and $%
\left\Vert u_{t_{1}}^{+}\right\Vert _{H_{k_{0}}}$. Indeed if we set $%
\left\Vert u_{t_{1}}\right\Vert =r\left( t_{1}\right) $ and $r\left(
t_{1}\right) $ satisfies the inequality 
\begin{equation*}
\widehat{a}\lambda _{k_{0}}+\widehat{b}\left( t_{1}\right) >\widehat{g}%
\left( t_{1}\right) r\left( t_{1}\right) ^{p_{1}\left( t_{1}\right) }>%
\widehat{a}\lambda _{k_{0}-1}+\widehat{b}\left( t_{1}\right) ,
\end{equation*}%
\ then each of solutions to the problem \ (\ref{5.1a}) -\ (\ref{5.1b})
satisfies one of the following statements:

1. If $u_{t_{1}}$\ lies in $H_{k_{0}}$ or in a small neighborhood of the
subspace $H_{k_{0}},$ then $\left\vert u_{k}\left( t\right) \right\vert
\nearrow \infty $ at $t\nearrow \infty $ for $k=\overline{1,k_{0}-1}$,
moreover since in this case $r\left( t\right) $ and $p_{1}\left( t\right) $
increase gradually and so in time $\widehat{g}\left( t\right) r\left(
t\right) ^{p_{1}\left( t\right) }$ is greater than each $\widehat{a}\lambda
_{k}+\widehat{b}\left( t\right) $ for $k\geq k_{0}$; i.e. in this case $%
\left\vert u_{k}\left( t\right) \right\vert $ gradually increases for all $k$%
;

2. If $u_{t_{1}}$\ lies in $H_{-k_{0}}$ or in a small neighborhood of the
subspace $H_{-k_{0}},$ then $\left\vert u_{k}\left( t\right) \right\vert
\searrow 0$ at $t\nearrow \infty $ for $k\geq k_{0}$, moreover since in this
case $r\left( t\right) $ decreases although $p_{1}\left( t\right) $
increases gradually so that in time $\widehat{g}\left( t\right) r\left(
t\right) ^{p_{1}\left( t\right) }$ is less than each $\widehat{a}\lambda
_{k}+\widehat{b}\left( t\right) $; i.e. in this case $\left\vert u_{k}\left(
t\right) \right\vert $ gradually decreases for all $k$;

3. If $u_{t_{1}}\in \boldsymbol{H}$ such that $\left\Vert
u_{t_{1}}^{-}\right\Vert _{H_{-k_{0}}}\approx \left\Vert
u_{t_{1}}^{+}\right\Vert _{H_{k_{0}}},$ then there is a relation between of $%
u_{t_{1}}^{-}\left( t\right) $ and $u_{t_{1}}^{+}$ such that the behavior of
the $u_{k}\left( t\right) $ is chaotic for all $k$ for which $u_{k}\left(
t_{1}\right) \neq 0$.

4. If $u_{t_{1}}\in \boldsymbol{H}$ such that $\left\Vert
u_{t_{1}}^{-}\right\Vert _{H_{-k_{0}}}$ and $\left\Vert
u_{t_{1}}^{+}\right\Vert _{H_{k_{0}}}$ are different, then the behavior of
solutions impossible to define, since its can be arbitrary due to vary of
the coefficients and the exponent.

Since the claims 1 and 2 were studied above, remains to study the claim 3
and the claim 4. We will use again of the representation of the formal
solutions (\ref{5.11}) to the problem\ (\ref{5.10}) 
\begin{equation*}
\left\vert u_{k}\left( t\right) \right\vert ^{2}=\exp \left\{ -2\left( 
\overset{t}{\underset{t_{1}}{\int }}\left( \widehat{a}\lambda _{k}+\widehat{b%
}\left( s\right) -\widehat{g}\left( s\right) r\left( s\right) ^{p_{1}\left(
s\right) }\right) ds\right) \right\} \left\vert u_{t_{1}k}\right\vert ^{2},
\end{equation*}%
$k=1,2,...$. It is known that $\left\vert u_{k}\left( t\right) \right\vert $
increases for $k:$ $1\leq k\leq k_{0}-1$ and decreases for $k\geq k_{0}$ by
virtue of Theorem \ref{Th_3}, depending on the difference $\widehat{a}%
\lambda _{k}+\widehat{b}\left( t_{1}\right) -\widehat{g}\left( t_{1}\right)
r\left( t_{1}\right) ^{p_{1}\left( t_{1}\right) }$. Since the $\widehat{g}%
\left( t\right) r\left( t\right) ^{p_{1}\left( t\right) }$ grows therefore
will be to take place $\widehat{a}\lambda _{k}+\widehat{b}\left( t\right) =%
\widehat{g}\left( t\right) r\left( t\right) ^{p_{1}\left( t\right) }$ for
some $k$ and $t>t_{2}$ in this case will appear the bifurcation that bring
to the appear of the chaos. Moreover as $r_{k_{0}}^{2}\left( t\right) \equiv
\left\Vert u_{k_{0}}^{-}\right\Vert _{H_{k_{0}}}^{2}\left( t\right)
+\left\Vert u_{k_{0}}^{+}\right\Vert _{H_{-k_{0}}}^{2}\left( t\right) $
holds near $t=t_{1}$, $\widehat{g}\left( t\right) r\left( t\right)
^{p_{1}\left( t\right) }$ can change depending on the behavior of $%
\left\vert u_{k}\left( t\right) \right\vert $, $\widehat{g}\left( t\right) $%
\ and $p_{1}\left( t\right) $ as $k\geq k_{0}$ and $k\leq k_{0}-1$ vary
consequently, the corresponding decomposation of $\boldsymbol{H}$ will
change, i.e. $\widehat{g}\left( t\right) r\left( t\right) ^{p_{1}\left(
t\right) }$ can become greater than $\widehat{a}\lambda _{k}+\widehat{b}%
\left( t\right) $ or less than $\widehat{a}\lambda _{k-1}+\widehat{b}\left(
t\right) $. Thus, the variation of $\rho \left( t\right) $ and $r\left(
t\right) $, and also their coupling are very complicated, as it depends on
relations among the behaviors of $u_{k}\left( t\right) $ in the case when $%
k\geq k_{0}$ and $k\leq k_{0}-1$, and also when change $k_{0}$, which may
give rise to chaos.

Now we will investigate the claim 4, i.e. the behavior of solutions of the
problem\ (\ref{5.1a}) - \ (\ref{5.1b}) when $u_{t_{1}}\in \boldsymbol{H}$
such that $\left\Vert u_{t_{1}k_{0}}^{-}\right\Vert _{H_{-k_{0}}}$ and $%
\left\Vert u_{t_{1}k_{0}}^{+}\right\Vert _{H_{k_{0}}}$ are different
moreover their relations in general not are such as in the claim 1 and claim
2. In the other words we will study properties of the resolving operator
corresponding to the problem (\ref{5.1a}) - \ (\ref{5.1b}). Toward this end,
we consider again to the following system of differential equations 
\begin{equation*}
\frac{1}{2}\frac{d}{dt}\left\langle u\left( t\right) ,w_{k}\right\rangle
+\left\langle \widehat{a}\nabla u\left( t\right) ,\nabla w_{k}\right\rangle +%
\widehat{b}\left( t\right) \left\langle u\left( t\right) ,w_{k}\right\rangle
-\widehat{g}\left( t\right) \left\Vert u\left( t\right) \right\Vert
^{p_{1}\left( t\right) }\left\langle u\left( t\right) ,w_{k}\right\rangle =0,
\end{equation*}%
where $\left\{ w_{k}\left( x\right) \right\} _{k=1}^{\infty }$ are
eigenfunctions of the Laplacian $-\Delta $ \ in $H_{0}^{1}\left( \Omega
\right) $ corresponding to the eigenvalues $\left\{ \lambda _{k}\right\}
_{k=1}^{\infty }$, respectively, by virtue of the imposed conditions.

Whence, it follows that 
\begin{equation}
\frac{1}{2}\frac{d}{dt}u_{k}\left( t\right) +\left( \widehat{a}\lambda _{k}+%
\widehat{b}\left( t\right) \right) u_{k}\left( t\right) -\widehat{g}\left(
t\right) \left\Vert u\left( t\right) \right\Vert ^{p_{1}\left( t\right)
}u_{k}\left( t\right) =0,\quad k=1,2,...  \label{5.12a}
\end{equation}

So, for the study of the behavior of solutions of the problem\ (\ref{5.1a})
- \ (\ref{5.1b}) under the condition $u_{t_{1}}\in
S_{r_{0}}^{H_{0}^{1}\left( \Omega \right) }\left( 0\right) $, $\widehat{a}%
\lambda _{k_{0}-1}+\widehat{b}\left( t\right) <\widehat{g}\left( t\right)
r_{0}^{p_{1}\left( t\right) }<\widehat{a}\lambda _{k_{0}}+\widehat{b}\left(
t\right) $ it is necessary to study the Cauchy problems in the case when $%
k=1,2,...,k_{0}-1$, and next the case $k\geq 1$. In the beginning we
consider the following Cauchy problem: 
\begin{equation}
\frac{1}{2}\frac{d}{dt}u_{k}\left( t\right) +\left( \widehat{a}\lambda _{k}+%
\widehat{b}\left( t\right) \right) u_{k}\left( t\right) -\widehat{g}\left(
t\right) \left\Vert u\left( t\right) \right\Vert ^{p_{1}\left( t\right)
}u_{k}\left( t\right) \ =0,  \label{5.13}
\end{equation}%
\begin{equation}
\left\langle u\left( t\right) ,w_{k}\right\rangle \left\vert \
_{t=t_{1}}\right. =u_{k}\left( t\right) \left\vert \ _{t=t_{1}}\right.
=u_{t_{1}k},\quad k=1,2,...,k_{0}-1  \label{5.14}
\end{equation}%
Whence, for some $t_{2}>t_{1}$ for $t\in \left[ t_{1},t_{2}\right) $ we have%
\begin{equation*}
\widehat{g}\left( t\right) \left\Vert u\left( t\right) \right\Vert
^{p_{1}\left( t\right) }\leq \widehat{g}\left( t\right) r_{0}^{p_{1}\left(
t\right) }+\varepsilon <\widehat{a}\lambda _{k_{0}}+\widehat{b}\left(
t\right) ,
\end{equation*}%
for some $\varepsilon >0$. Indeed, we have from \ (\ref{5.2}) - (\ref{5.3})
that 
\begin{eqnarray*}
\frac{d}{dt}u_{k}\left( t\right) +2\left[ (\widehat{a}\lambda _{k}+\widehat{b%
}\left( t\right) -\widehat{g}\left( t\right) \left\Vert u\left( t\right)
\right\Vert ^{p_{1}\left( t\right) }\right] u_{k}\left( t\right) &=&0, \\
u_{k}\left( t_{1}\right) &=&u_{t_{1}k},
\end{eqnarray*}%
which leads to the formal solution of the Cauchy problem 
\begin{equation}
u_{k}\left( t\right) =\exp \left\{ -2\underset{t_{1}}{\overset{t}{\int }}(%
\widehat{a}\lambda _{k}+\widehat{b}\left( s\right) -\widehat{g}\left(
s\right) r\left( s\right) ^{p_{1}\left( s\right) })ds\right\} u_{t_{1}k}.
\label{5.15}
\end{equation}%
Consequently, if $\widehat{a}\lambda _{k_{0}-1}+\widehat{b}\left(
t_{1}\right) <\widehat{g}\left( t\right) r_{t_{1}}^{p_{1}\left( t\right) }$,
then $u_{k}\left( t,x\right) $ increases in the vicinity of $%
u_{t_{1}k}\left( x\right) $ if $u_{t_{1}k}\left( x\right) >0$ for $%
k=1,2,...,k_{0}-1$ and the part $\left\Vert P_{k_{0}}u_{t_{1}}\right\Vert $
is enough greater than $\left\Vert Q_{k_{0}}u_{t_{1}}\right\Vert $.

Let $u_{t_{1}}\in B_{r_{0}}^{H_{0}^{1}}\left( 0\right) $ and $\left\Vert
u_{t_{1}}\right\Vert =r_{0}$, then the above expression implies that the
behavior of the solution $u_{k}\left( t\right) $ depends on the relation
between $\left\Vert P_{k}u_{t_{1}}\right\Vert $ and $\left\Vert
Q_{k}u_{t_{1}}\right\Vert $, and also depends of the relationship between $%
\widehat{g}\left( t_{1}\right) r_{0}^{p_{1}\left( t_{1}\right) }$ and $%
\widehat{a}\lambda _{k}+\widehat{b}\left( t_{1}\right) $.

Now for study of the behavior of the $u\left( t,x\right) =\underset{k=1}{%
\overset{\infty }{\sum }}u_{k}\left( t\right) w_{k}\left( x\right) $
consider the system of equations that one can write as

\begin{equation*}
\frac{d}{dt}u_{k}\left( t\right) +2\left( \widehat{a}\lambda _{k_{0}-1}+%
\widehat{b}\left( t\right) -\widehat{g}\left( t\right) r_{0}^{p_{1}\left(
t\right) }\right) u_{k}\left( t\right) ,\ 
\end{equation*}%
\begin{equation}
u_{k}\left( t_{1}\right) =u_{t_{1}k},\quad \ k=1,2,...  \label{5.16}
\end{equation}%
according to assumption.

Let there exists $\lambda _{k_{0}}$, $k\geq 2$ such that the inequalities $%
\widehat{a}\lambda _{k_{0}-1}+\widehat{b}\left( t_{1}\right) <\widehat{g}%
\left( t_{1}\right) r_{0}^{p_{1}\left( t_{1}\right) }<\widehat{a}\lambda
_{k_{0}}+\widehat{b}\left( t_{1}\right) $ hold, whence isn't difficult to
see that this system of equations is of interest only for the cases $k\geq
k_{0}$, $k\leq k_{0}-1$ and $k=k_{0}$ separately. If $k\geq k_{0}$, this
part of the system has a solution as $t\leq t_{2}\left( k,r_{0},p_{1}\right) 
$ for some $t_{2}\left( k,r_{0},p_{1}\right) >t_{1}$ if $\left\Vert
Q_{k_{0}}u_{t_{1}}\right\Vert $ is sufficiently greater than $\left\Vert
P_{k_{0}}u_{t_{1}}\right\Vert $. Formally we can determine the solution of
each equation from \ (\ref{5.9})-(\ref{5.9a}) in the following form: 
\begin{equation}
u_{k}\left( t\right) =\exp \left\{ -2\left( \overset{t}{\underset{t_{1}}{%
\int }}\widehat{a}\lambda _{k}+\widehat{b}\left( s\right) -\widehat{g}\left(
s\right) r_{0}^{p_{1}\left( s\right) }ds\right) \right\} u_{t_{1}k}\ .
\label{5.17}
\end{equation}%
Thus, considering the problem (\ref{5.9})-(\ref{5.9a}) for $k\leq k_{0}-1$, (%
\ref{5.12a}) implies that 
\begin{equation*}
\left\vert u_{k}\left( t\right) \right\vert =\exp \left\{ -2\left( \overset{t%
}{\underset{t_{1}}{\int }}\widehat{a}\lambda _{k}+\widehat{b}\left( s\right)
-\widehat{g}\left( s\right) r_{0}^{p_{1}\left( s\right) }ds\right) \right\}
\left\vert u_{t_{1}k}\right\vert \geq \quad
\end{equation*}%
\begin{equation*}
\geq \exp \left\{ 2\left( \widehat{g}\left( t_{1}\right) r_{0}^{p_{1}\left(
t_{1}\right) }-\widehat{a}\lambda _{k}+\widehat{b}\left( t_{1}\right)
\right) t\right\} \left\vert u_{t_{1}k}\right\vert ,
\end{equation*}%
as $\widehat{a}\lambda _{k}+\widehat{b}\left( t\right) <\widehat{g}\left(
t\right) r_{0}^{p_{1}\left( t\right) }\ $ for$\ 1\leq k\leq k_{0}-1$ and
some $t>t_{1}.$ Consequently, the sequence $\left\vert u_{k}\left( t\right)
\right\vert $ increases for each $k:$ $1\leq k\leq k_{0}-1$, if the part $%
\left\Vert P_{k_{0}}u_{t_{1}}\right\Vert $ is sufficiently greater than $%
\left\Vert Q_{k_{0}}u_{t_{1}}\right\Vert $, and that renders the increase of 
$r\left( t\right) $.

Consider the representation of the formal solutions \ (\ref{5.16}) to the
problem\ (\ref{5.9})-(\ref{5.9a}): 
\begin{equation*}
u_{k}\left( t\right) =\exp \left\{ -2\left( \overset{t}{\underset{t_{1}}{%
\int }}\widehat{a}\lambda _{k}+\widehat{b}\left( s\right) -\widehat{g}\left(
s\right) r_{0}^{p_{1}\left( s\right) }ds\right) \right\} u_{t_{1}k},\quad
k=1,2,....
\end{equation*}%
It is known that $\left\vert u_{k}\left( t\right) \right\vert $ increases
for $k:$ $1\leq k\leq k_{0}-1$ and decreases for $k\geq k_{0}$ by virtue of
Theorem \ref{Th_3}, depending on the difference $\widehat{a}\lambda _{k}+%
\widehat{b}\left( t_{1}\right) -\widehat{g}\left( t_{1}\right)
r_{0}^{p_{1}\left( t_{1}\right) }$ and the relation between $\left\Vert
P_{k_{0}}u_{0}\right\Vert $ and $\left\Vert Q_{k_{0}}u_{0}\right\Vert $.
Thus, in order to investigate of the system (\ref{5.16}) for each $k$ we
need to study the expression $\widehat{a}\lambda _{k}+\widehat{b}\left(
t_{1}\right) -\widehat{g}\left( t_{1}\right) r_{0}^{p_{1}\left( t_{1}\right)
}$ which can be negative or positive due to the conditions imposed. But the
behavior of functions $|u_{k}(t)|$ cannot exactly explain the behavior of
functions $u_{k}(t)$, and also the behavior of the solution $u(t,x)$.
Consequently, we need to study the behavior of functions $u_{k}(t)$ in
greater detail.

So, from the representation (\ref{5.17}) under the corresponding relation
between $\left\Vert P_{k_{0}}u_{t_{1}}\right\Vert $ and $\left\Vert
Q_{k_{0}}u_{t_{1}}\right\Vert $ follows that if $u_{t_{1}k}\geq 0$ then $%
u_{k}\left( t\right) \geq 0$ and in addition if $\widehat{a}\lambda _{k}+%
\widehat{b}\left( t_{1}\right) -\widehat{g}\left( t_{1}\right)
r_{0}^{p_{1}\left( t_{1}\right) }>0$ then one can see that $u_{k}\left(
t\right) $ will decreases at least in some vicinity of $u_{t_{1}k}$, and if $%
\widehat{a}\lambda _{k}+\widehat{b}\left( t_{1}\right) -\widehat{g}\left(
t_{1}\right) r_{0}^{p_{1}\left( t_{1}\right) }<0$, then when $u_{t_{1}k}\geq
0$ one can see that $u_{k}\left( t\right) $ increases at least in a vicinity
of $u_{t_{1}k}$. But if $\widehat{a}\lambda _{k}+\widehat{b}\left(
t_{1}\right) -\widehat{g}\left( t_{1}\right) r_{0}^{p_{1}\left( t_{1}\right)
}=0$ then not is possible to understand how will vary the $u_{k}\left(
t\right) $ in vicinity of $u_{t_{1}k}$.

Thus, we need to investigate the behavior of $r\left( t\right) $ for various
function $u_{t_{1}}\left( x\right) $ in the case when $\left\Vert
u_{t_{1}}\right\Vert =r_{0}$, and more exactly for various function $%
u_{0}\left( x\right) $ since $u_{t_{1}}\left( x\right) $ depends on $%
u_{0}\left( x\right) $. Consequently, the behavior of $r\left( t\right) $
essentially depends on the geometrical selections of $u_{k}\left(
t_{1}\right) $ (more exactly on the geometrical selections of initial
function $u_{0}$).

In the above assumptions the solution of examined problem (\ref{5.1a})-(\ref%
{5.1b}) and of the initial function at the moment $t_{1}$ that is obtained
from initial data $u_{0}$ we can represent with the expressions $u\left(
t,x\right) =\underset{k\geq 1}{\sum }u_{k}\left( t\right) w_{k}\left(
x\right) $ and $u_{t_{1}}\left( x\right) =\underset{k\geq 1}{\sum }%
u_{t_{1}k}w_{k}\left( x\right) $, respectively. The provided analysis shows
that for the study of the behavior of solutions of the problem in detail it
is necessary to use the above expressions.

So, let $r_{0}>0$, then there exists $k_{0}\geq 2$ such that 
\begin{equation*}
\widehat{a}\lambda _{k_{0}-1}+\widehat{b}\left( t_{1}\right) <\widehat{g}%
\left( t_{1}\right) r_{0}^{p_{1}\left( t_{1}\right) }\leq \widehat{a}\lambda
_{k_{0}}+\widehat{b}\left( t_{1}\right)
\end{equation*}%
and $\left\Vert u_{t_{1}}\right\Vert =r_{0}$.

Then by use the orthogonal splitting $u\left( t\right) =P_{k_{0}}u\left(
t\right) +Q_{k_{0}}u\left( t\right) $, we obtain the following
represetations: 
\begin{equation*}
u_{t_{1}}\left( x\right) =\underset{k\geq 1}{\sum }u_{t_{1}k}w_{k}\left(
x\right) =\underset{k_{0}-1\geq k\geq 1}{\sum }u_{t_{1}k}w_{k}\left(
x\right) +\underset{k\geq k_{0}}{\sum }u_{t_{1}k}w_{k}\left( x\right)
\end{equation*}%
and 
\begin{equation}
u\left( t,x\right) =\underset{k\geq 1}{\sum }u_{k}\left( t\right)
w_{k}\left( x\right) =\underset{k_{0}-1\geq k\geq 1}{\sum }u_{k}\left(
t\right) w_{k}\left( x\right) +\underset{k\geq k_{0}}{\sum }u_{k}\left(
t\right) w_{k}\left( x\right) .  \label{5.18}
\end{equation}

From the relations between of the expressions $\widehat{g}\left( t\right)
r\left( t\right) ^{p_{1}\left( t\right) }$and $\widehat{a}\lambda _{k_{0}-1}+%
\widehat{b}\left( t\right) $ follows there exist $t_{2}>t_{1}$ and $%
t_{3}>t_{1}$ such that $P_{k_{0}}u\left( t\right) $ can increases in $\left(
t_{1},t_{2}\right) $ and $Q_{k_{0}}u\left( t\right) $ decreases in $\left(
t_{1},t_{3}\right) $ in the formula (\ref{5.18}) by virtue of the (\ref{5.17}%
). Then if $\min \left\{ t_{2},t_{3}\right\} =t_{2}$, then for $t>$ $t_{2}$
these summands can behave quite differently. Here there exist the following
possibilities: 1) the functional $\left\Vert u\right\Vert \left( t\right) $
becomes greater than $r_{0}$ for $t\geq t_{2}$, moreover $\widehat{g}\left(
t\right) r\left( t\right) ^{p_{1}\left( t\right) }\geq \widehat{a}\lambda
_{k_{0}}+\widehat{b}\left( t\right) $ for $t>t_{3}$, so the orthogonal
splitting $u=P_{k_{0}}u+Q_{k_{0}}u$ changes and becomes, at least, $%
u=P_{k_{0}+1}u+Q_{k_{0}+1}u$; 2) $Q_{k_{0}}u$ decreases up to a point where $%
\left\Vert u\right\Vert \left( t\right) $ is smaller than $r_{0}$ for $t\geq
t_{2}$, moreover $\widehat{g}\left( t\right) r\left( t\right) ^{p_{1}\left(
t\right) }\leq \widehat{a}\lambda _{k_{0}-1}+\widehat{b}\left( t\right) $
for $t>t_{2}$, so the orthogonal splitting $u=P_{k_{0}}u+Q_{k_{0}}u$ changes
and becomes, at least, $u=P_{k_{0}-1}u+Q_{k_{0}-1}u$; 3) there exist a $%
t_{4}\geq t_{2}$ and an $R_{0}\geq \left\Vert P_{k_{0}}u\right\Vert \geq
R_{1}>0$ such that beginning at $t_{4}$ the changes of $P_{k_{0}}u$ and $%
Q_{k_{0}}u$\ become such that 
\begin{equation}
r^{2}\left( t\right) =\left\Vert u\left( t\right) \right\Vert ^{2}=\underset{%
k_{0}-1\geq k\geq 1}{\sum }\left\vert u_{k}\left( t\right) \right\vert ^{2}+%
\underset{k\geq k_{0}}{\sum }\left\vert u_{k}\left( t\right) \right\vert ^{2}
\label{5.19}
\end{equation}%
satisfies $R_{1}\leq r\left( t\right) \leq R_{0}$ for $t\geq t_{4}$.

Consider case 1 where we have the following possibilities: \textit{(i)}\ $%
P_{k_{0}}u$ increases with such velocity that $\left\Vert u\right\Vert
\left( t\right) \nearrow \infty $, which can occur when $u_{t_{1}}\left(
x\right) $ is chosen near the subspace $H_{k_{0}}$ (this scenario is studied
in Theorem \ref{Th_3}); \textit{(ii) }the\textit{\ }rate of growth of $%
P_{k_{0}}u$ diminishes beginning at time $t$ and the function $u\left(
t,x\right) $ behaves as in case 3, which we will explain in what follows.
Case 2 has 2 variants: \textit{(iii)}\ $Q_{k_{0}}u$ decreases with such
velocity that $\left\Vert u\right\Vert \left( t\right) \searrow 0$ leading
to the inequality $\widehat{g}\left( t\right) r\left( t\right) ^{p_{1}\left(
t\right) }<\widehat{a}\lambda _{1}+\widehat{b}\left( t\right) $, which can
take place when $u_{t_{1}}\left( x\right) $ is chosen in near the subspace $%
H_{-k_{0}}$ (this variant is also studied in Theorem \ref{Th_3}); \textit{%
(iv) }the\textit{\ }rate of decrease of $Q_{k_{0}}u$ diminishes beginning at
some time $t$ and leads to case 1\textit{(ii)}.

We should be noted all of the above cases depends of the geometrical
location of the initial function on the sphere $S_{r_{0}}^{H_{0}^{1}}\left(
0\right) \subset H_{0}^{1}\left( \Omega \right) $.

Consequently, it remains only to investigate case 3. Therefore, we can
choose special initial data $u\left( 0\right) =u_{0}$ such that into the
time $t_{1}$ the functions $u\left( t_{1}\right) $ be on the sphere $%
S_{r_{0}}^{H_{0}^{1}}\left( 0\right) $ and try to explain case 3 for such
functions. So, let $P_{k_{0}}u_{t_{1}}=u_{t_{1}k_{0}-1}w_{k_{0}-1}$, i.e. $%
u_{t_{1}}\left( x\right) =u_{t_{1}k_{0}-1}w_{k_{0}-1}\left( x\right)
+Q_{k_{0}}u_{t_{1}}\left( x\right) $ and $\left\Vert u_{t_{1}}\right\Vert
=r_{0}$. Then we obtain the change of the function $u\left( t,x\right) $
happen by the following way: $u_{k_{0}-1}\left( t\right) $ changes so that $%
\left\vert u_{k_{0}-1}\left( t\right) \right\vert $ increases with $t$ and $%
\widehat{g}\left( t\right) \left\vert u_{k_{0}-1}\left( t\right) \right\vert
^{p_{1}\left( t\right) }\longrightarrow \widehat{a}\lambda _{k_{0}-1}+%
\widehat{b}\left( t\right) $ when $t\nearrow \infty $, and $Q_{k_{0}}u\left(
t\right) $ changes so that $\left\Vert Q_{k_{0}}u\left( t\right) \right\Vert 
$ decreases although $p_{1}\left( t\right) $ increases with increasing $t$
and therefore $\left\Vert Q_{k_{0}}u\left( t\right) \right\Vert
\longrightarrow 0$ when $t\nearrow \infty $. Hence, $\widehat{g}\left(
t\right) \left\Vert u\left( t\right) \right\Vert ^{p_{1}\left( t\right)
}\searrow \widehat{a}\lambda _{k_{0}-1}+\widehat{b}\left( t\right) $ as $%
t\nearrow \infty $. In other words, the increase of $\left\Vert
P_{k_{0}}u\right\Vert $ and decrease of $\left\Vert Q_{k_{0}}u\left(
t\right) \right\Vert $ compensate for each other in such a way that this
process leads to the case described above.

We should be noted in order to the functions $u_{t_{1}}$ satisfy of these
conditions is necessary the initial data $u_{0}$ to choose by corresponding
way, namely using of the representation (\ref{4.9b}), that isn't difficult
see. Obviously that the function $u_{0}$ must be has also of same
representation.

Thus, it isn't difficult to see in order to obtain the above result, we need
to select $u_{t_{1}k_{0}-1}$ (and also $u_{0k_{0}-1}$) in near of $%
H_{-k_{0}} $, which depends on the given%
\begin{equation}
r_{0}:\widehat{a}\lambda _{k_{0}-1}+\widehat{b}\left( t_{1}\right) <\widehat{%
g}\left( t_{1}\right) r_{0}^{p_{1}\left( t_{1}\right) }\leq \widehat{a}%
\lambda _{k_{0}}+\widehat{b}\left( t_{1}\right) .  \label{5.20a}
\end{equation}%
Accordingly it follows in the case when $P_{k_{0}}u\left( t\right) $
increases, the corresponding $u_{t_{1}k}$, $1\leq k\leq $ $k_{0}-1$, must be
chosen as done previously. Moreover, there is a $\lambda _{j_{0}}$ such that 
$\widehat{g}\left( t\right) \left\Vert P_{k_{0}}u\left( t\right) \right\Vert
^{p_{1}\left( t\right) }\nearrow \widehat{a}\lambda _{j_{0}}+\widehat{b}%
\left( t\right) $ when $t\nearrow \infty $, where 
\begin{equation*}
\lambda _{j_{0}}=\inf \left\{ \lambda _{k}\left\vert \ 1\leq k\leq
k_{0}-1,\right. u_{t_{1}k}\neq 0\right\} .
\end{equation*}

Consequently, we arrive that there exists a "double cone" with the "vertex
at zero" that contains the subspace $H_{-k_{0}}$ and all elements are
contained in some neighborhood of $H_{-k_{0}}$, where the maximal distance
between of the elements of this subset and the subspace $H_{-k_{0}}$ depends
on the given $r_{0}$, and also on the coeffitiens and the exponent from the
inequality (\ref{5.20a}). If we denote this subset by $\widetilde{H}\subset
H_{0}^{1}\left( \Omega \right) $, then we have that any subset of $%
\widetilde{H}\cap \left\{ B_{r_{1}}^{H_{0}^{1}}\left( 0\right)
-B_{r_{2}}^{H_{0}^{1}}\left( 0\right) ,\ r_{1}>r_{2}>0\right\} $ converges
to a set, which we can define as $H_{k_{0}}\cap
B_{r_{j_{0}}}^{H_{0}^{1}}\left( 0\right) $, where $r_{j_{0}}=\widehat{a}%
\lambda _{j_{0}}+\widehat{b}\left( t_{1}\right) $ and$\ r_{1},r_{2}$ are
some numbers with 
\begin{equation*}
\widehat{a}\lambda _{k_{1}-1}+\widehat{b}\left( t_{1}\right) <\widehat{g}%
\left( t_{1}\right) r_{2}^{p_{1}\left( t_{1}\right) }\leq \widehat{a}\lambda
_{k_{1}}+\widehat{b}\left( t_{1}\right)
\end{equation*}%
and there is a $\lambda _{j_{1}}=\inf \left\{ \lambda _{k}\left\vert \ 1\leq
k\leq k_{1}-1,\right. u_{t_{1}k}\neq 0\right\} $ and $k_{1}=k_{1}\left(
r_{1}\right) $. This shows that $H_{k_{0}}\cap
B_{r_{j_{1}}}^{H_{0}^{1}}\left( 0\right) $ is a subset of a
finite-dimensional space and it is local attractor in some sense, where $%
r_{j_{1}}=\widehat{a}\lambda _{j_{1}}+\widehat{b}\left( t_{1}\right) $.

Thus, we have proved the following result.

\begin{theorem}
\label{Th_5}Let all the above conditions hold and $u_{t_{1}}\in
H_{0}^{1}\left( \Omega \right) $ satisfies the inequality $\widehat{a}%
\lambda _{k_{0}-1}+\widehat{b}\left( t_{1}\right) <\widehat{g}\left(
t_{1}\right) \left\Vert u\left( t_{1}\right) \right\Vert
_{H_{0}^{1}}^{p_{1}\left( t_{1}\right) }\leq \widehat{a}\lambda _{k_{0}}+%
\widehat{b}\left( t_{1}\right) $.\footnote{%
The relation between $u_{t_{1}}$ and $u_{0}$ is given by formula (\ref{4.9b})%
} Then each solution to the problem \ (\ref{5.1a}) -\ (\ref{5.1b}) satisfies
one of the following properties:

1. If $u_{_{t_{1}}}$\ lies in $H_{k_{0}}$ or in a sufficiently small
neighborhood of the subspace $H_{k_{0}},$ then $\left\vert u_{k}\left(
t\right) \right\vert \nearrow \infty $ as $t\nearrow \infty $ for $k=%
\overline{1,k_{0}-1}$; moreover, in this case $r\left( t\right) =\left\Vert
u\left( t\right) \right\Vert $ increases and $\widehat{g}\left( t\right)
r\left( t\right) ^{p_{1}\left( t\right) }$ gradually becomes greater than $%
\widehat{a}\lambda _{k}+\widehat{b}\left( t\right) $ for each $\lambda _{k}:$
$k\geq k_{0}$, $\left\vert u_{k}\left( t\right) \right\vert $ gradually
increases for all $k$;

2. If $u_{_{t_{1}}}$\ lies in $H_{-k_{0}}$ or in a small neighborhood of the
subspace $H_{-k_{0}},$ then $\left\vert u_{k}\left( t\right) \right\vert
\searrow 0$ as $t\nearrow \infty $ for $k\geq k_{0}$; moreover, since in
this case $r\left( t\right) $ decreases and although $p_{1}\left( t\right) $
increases but $\widehat{g}\left( t\right) r\left( t\right) ^{p_{1}\left(
t\right) }$ gradually becomes less than $\widehat{a}\lambda _{k}+\widehat{b}%
\left( t\right) $ for each $\lambda _{k}$, $\left\vert u_{k}\left( t\right)
\right\vert $ gradually decreases for all $k$;

3. If $u_{_{t_{1}}}\in \boldsymbol{H}$, $\left\Vert
P_{k_{0}}u_{_{t_{1}}}\right\Vert \ll \left\Vert
Q_{k_{0}}u_{_{t_{1}}}\right\Vert $ and if there are small numbers%
\begin{equation*}
\delta \left( \widehat{a}\lambda _{k_{0}}+\widehat{b}\left( t_{1}\right)
\right) >\epsilon \left( \widehat{a}\lambda _{k_{0}}+\widehat{b}\left(
t_{1}\right) \right) >0
\end{equation*}%
such that the Hausdorff distance satisfies 
\begin{equation}
\epsilon \leq d\left( H_{-k_{0}};\left\{ u_{t_{1}k}^{+}\left\vert \ k=%
\overline{1,k_{0}-1}\right. \right\} \right) \leq \delta ,  \label{5.20}
\end{equation}%
then $u\left( t,x\right) $ is chaotic for sufficient large $t$. In addition,
if $\left\Vert P_{k_{0}}u_{t_{1}}\right\Vert _{H_{-k_{0}}}\approx \left\Vert
Q_{k_{0}}u_{t_{1}}\right\Vert _{H_{k_{0}}},$ then there is a relationship
between $u_{t_{1}k_{0}}^{-}\left( t\right) $ and $u_{t_{1}k_{0}}^{+}$ for
which the behavior of the $u_{k}\left( t\right) $ is chaotic for all $k$
satisfying $u_{k}\left( 0\right) \neq 0$.
\end{theorem}

\begin{remark}
If property 3 of the above theorem obtains, then the following claim is
reasonable: for any $\lambda _{k_{0}}$ there is a subset $%
D_{r_{k_{0}}}\subset H_{0}^{1}\left( \Omega \right) $ for which (\ref{5.20})
holds and for any $u_{t_{1}}\in D_{r_{k_{0}}}$ the corresponding solution $%
u\left( t\right) $ satisfies the condition 
\begin{equation*}
\widehat{a}\lambda _{j_{0}}+\widehat{b}\left( t_{1}\right) \leq \widehat{g}%
\left( t\right) r\left( t\right) ^{p_{1}\left( t\right) }<\widehat{a}\lambda
_{k_{0}}+\widehat{b}\left( t_{1}\right) \lambda _{k_{0}}
\end{equation*}%
for any $t>0$, where $r_{k_{0}}=\widehat{a}\lambda _{k_{0}}+\widehat{b}%
\left( t_{1}\right) $ and \ 
\begin{equation*}
\lim {}_{t\nearrow \infty }\ \widehat{g}\left( t\right) r\left( t\right)
^{p_{1}\left( t\right) }=\widehat{a}\lambda _{j_{0}}+\widehat{b}\left(
t_{1}\right) ,\quad u_{0k}\neq 0
\end{equation*}%
\ then there is an absorbing chaotic set in $L^{2}\left( \Omega \right) $,
where 
\begin{equation*}
\lambda _{j_{0}}=\inf \left\{ \lambda _{k}\left\vert \ 1\leq k\leq
k_{0}-1\right. \right\} .
\end{equation*}
\end{remark}

\subsection{\label{Con_4}Conclusion on the dynamics of cancer}

The obtained results show that the propagation of cancer essentially depends
on the initial value, i.e. what moves the start of the destruction
(degeneration) and of which initial mass covers. It is naturally that the
dynamics of propagation of cancer essential depends on the immune systems of
the body. Naturally, the dynamics of propagation of cancer essential depends
on the immune systems of the body, consequently on its powers to act and to
stop the degeneration process, and also the sufficient stability of its act.
So, we can classify the initial date of the space of initial dates beginning
with their role in participation in the process of the propagation of
cancer. Thus one can determine 3 classes of the initial dates that have a
principal role in the studied process. The first class contains such initial
dates, under which the process has the successful result. (The cases
investigated in sections 2.1 and 2.3 contain this class.) The second class
contains such initial dates, under which the process has the worse result.
(The case investigated in section 3.1 contain this class.) The third class
is found between classes 1 and 2 and contains such initial dates, under
which the process of the dynamics of propagation of cancer is sufficiently
complicated. The appearance of chaos in the dynamics at the process of the
propagation of cancer is connected with the initial date and stationary part
of the problem. Since at the initial value from this class due to the
stationary part of the examined problem arise the bifurcations, which turn
into for the examined process to the cascade of bifurcations, therefore
exist such absorbing manifold, that the associated dynamics become chaotic.
Consequently, in this case, the long-time behavior of the trajectory of the
dynamics of propagation of cancer can be mainly chaotic, therefore it become
undefinable.

\section{Open Problems}

The above discuss shows: in order to the mathematical model represented in
this article could describe of the dynamics of the cancer more exactly is
necessary that all functions in the equation (\ref{1.1}) depends from the
both variables $t$ and $x$.

\begin{problem}
Do possible to study the long-time behavior of solutions of the problem (\ref%
{1.1})-(\ref{1.3}) in the above-mentioned conditions?
\end{problem}

\begin{problem}
Do possible to study the control problem for the dynamics of the cancer
using the represented in this article mathematical model: the problem (\ref%
{1.1})-(\ref{1.3})?

We should be noted the control function must depends from the both variables
and have the variable exponent, and also it is need to determine of the
corresponding class of the control functions.
\end{problem}

\bigskip

\bigskip

\end{document}